\documentclass[twocolumn]{autart}    

\usepackage{graphicx}
\usepackage[T1]{fontenc}
\usepackage{amsmath}
\usepackage{amssymb}
\usepackage{amsfonts}
\usepackage{color}
\usepackage{float}
\usepackage{subfigure}
\usepackage{url}
\usepackage{calligra}
\usepackage{adjustbox}

\newtheorem{definition}{Definition}

\newtheorem{theorem}{Theorem}

\begin{document}

\begin{frontmatter}

\title{Distributed Event-Triggered Nash Equilibrium Seeking\\ for Noncooperative Games} 


\thanks[footnoteinfo]{Corresponding author: V.~H.~P.~Rodrigues.}

\author[first]{Victor Hugo Pereira Rodrigues\thanksref{footnoteinfo}}\ead{rodrigues.vhp@gmail.com}\,,
\author[first]{Tiago Roux Oliveira}\ead{tiagoroux@uerj.br}\,,\linebreak 
\author[second]{Miroslav Krsti{\' c}}\ead{mkrstic@ucsd.edu}\,, and 
\author[third]{Tamer Ba{\c s}ar}\ead{basar1@illinois.edu}\,. 
       
\address[first]{State University of Rio de Janeiro (UERJ), Rio de Janeiro -- RJ, Brazil}             
\address[second]{University of California San Diego (UCSD), San Diego -- CA, USA}        
\address[third]{University of Illinois Urbana-Champaign, Urbana -- IL, USA}
          

\begin{abstract}
We propose locally convergent Nash equilibrium seeking algorithms for $N$-player noncooperative games, which use distributed event-triggered pseudo-gradient estimates. The proposed approach employs sinusoidal perturbations to estimate the pseudo-gradients of unknown quadratic payoff functions. This is the first instance of noncooperative games being tackled in a model-free fashion with event-triggered extremum seeking. Each player evaluates independently the deviation between the corresponding current pseudo-gradient estimate and its last broadcasted value
from the event-triggering mechanism to tune individually the player action, while they preserve collectively the closed-loop stability/convergence. 
%
We guarantee Zeno behavior avoidance by establishing a minimum dwell-time to avoid infinitely fast switching. In particular, the stability analysis is carried out using Lyapunov's method and averaging for systems with discontinuous right-hand sides. We quantify the size of the ultimate small residual sets around the Nash equilibrium and illustrate the theoretical results numerically on an oligopoly setting.
\end{abstract}

\end{frontmatter}

\section{Introduction}

\textcolor{black}{
Game theory provides a theoretical framework for studying social situations among competing players and using mathematical models of strategic interaction among rational decision-makers \cite{Tirole:1991,BZ:2018}.
Game-theoretic approaches to designing, modeling, and optimizing emerging engineering systems, biological
behaviors and mathematical finance make this topic of research an extremely important one in many fields with a wide range of applications \cite{Basar:2019,Sastry:2013,BZ2:2018}. Indeed, there exist numerous studies of both theoretical and practical nature in the literature on differential games
\cite{SH:1969,PG:1988,SB:1991,WSBX:2018,CZ:2018,AS:2019,ATM:2020}. }

\textcolor{black}{
In this context, we can view games roughly in two categories: \textit{cooperative} and \textit{noncooperative} games \cite{Basar:1999}. A game is cooperative if the players are able to form binding commitments externally enforced (\textit{e.g.}, through contract law), resulting in collective payoffs. A game is noncooperative if players cannot form alliances or if all agreements need to be self-enforcing (\textit{e.g.}, through credible threats), focusing on predicting individual players' actions and payoffs and analyzing the \textit{Nash equilibrium} \cite{N:1951}. Nash equilibrium is an outcome which, once achieved or reached, means no player can improve upon its payoff by modifying its decision unilaterally \cite{Basar:1999}. }

\textcolor{black}{
The development of algorithms to achieve convergence
to a Nash equilibrium has been a focus of researchers for
several decades \cite{LB:1987,B:1987}. Some papers have also looked at learning aspects of various update schemes for reaching Nash equilibrium \cite{ZTB:2013}. Related to extremum seeking (ES), the authors in \cite{FKB:2012} study the problem of computing, in real time, the Nash equilibrium of static noncooperative games with $N$ players by employing a non-model based approach. By utilizing ES \cite{KW:2000} with sinusoidal perturbations, the players achieve stable, local attainment of their Nash strategies without the need for any model information. Thereafter, the expansion of the algorithm proposed in \cite{FKB:2012} to different domains was given for delayed and PDE systems \cite{JOTA:2021,CSm:2024} as well as fixed-time convergence \cite{Poveda:2023,VSS:2024}.
}

\textcolor{black}{
On the other hand, in the current technological age of network science, researchers are focusing on decreasing costs by designing  fast and reliable communication schemes where the plant and controller might not be physically connected, or might even be in different geographical locations. These networked control systems  offer advantages in the financial cost of installation and maintenance \cite{ZHGDDYP:2020}. However, one of their major disadvantages is the resulting high-traffic congestion, which can lead to transmission delays and  packet dropouts, {\it i.e.}, data may be lost while in transit through the network \cite{HNX:2007}. These issues are highly related to limited resource or available communication channels' bandwidth. To alleviate or mitigate this problem, Event-Triggered Controllers (ETC) can be used.}

\textcolor{black}{
ETC executes the control task, non-periodically, in response to a triggering condition designed as a function of the plant's state. Besides the asymptotic stability properties \cite{T:2007}, this strategy reduces control effort since the control update and data communication only occur when the error between the current state and the equilibrium set exceeds a certain value that might induce instability \cite{BH:2013}. Pioneering works towards the development of resource-aware control design include the construction of digital computer design \cite{s9b}, the event-based PID design  \cite{s9}  and the event-based controller for stochastic systems \cite{s8}. Works dedicated to extension of event-based control for   networked control systems with a high level of complexity exist  for both linear \cite{s1,HJT:2012,s5} and nonlinear systems \cite{T:2007,APDN:2016}. 
Among others, results on event-based control  deal with the robustness against the effect of  possible  perturbations \cite{s13,s14} or parametric uncertainties \cite{s16}. In  \cite{ZLJ:2021},  ETC is designed to satisfy a cyclic-small-gain condition such that the stabilization of a class of nonlinear time-delay systems is guaranteed. Further, the authors in \cite{CL:2019} propose  distributed event-triggered leaderless and leader-follower consensus control laws for general linear multi-agent networks. An event-triggered output-feedback  design \cite{APDN:2016} has been employed aiming to stabilize a class of nonlinear systems by combining techniques from event-triggered and time-triggered control. Research in this area spans various control-estimation designs and system complexities, addressing robustness and stability concerns \cite{IB:2005,IB:2010a,IB:2010,APDN:2016,ZLJ:2021,CL:2019}.
}

\textcolor{black}{
Many practical engineering problems that arise in industry, involving networking in areas such as network virtualization, software defined networks, cloud computing, the Internet of Things, context-aware networks, green communications, and security  can be modeled by a game-theoretic approach with
distributively connectivity through networks by using their resources \cite{Basar:2019,Sastry:2013,AB:2011}. Hence, the motivation for employing real-time optimization to improve such engineering processes commonly modeled by this combination of game theory and event-triggered architectures in networked-based framework should be clear and highly demanding. There even exist publications on these topics \cite{Mazo_Tabuada,Wang_Lemmon,Johansson,Xinghuo_Yu}, but the literature has not yet addressed in this context extremum seeking feedback \cite{KW:2000}.}

\textcolor{black}{
In fact, it appears to be quite challenging to address this problem. Despite the large number of publications on ETC mentioned in the previous paragraphs and also the consolidation of ES results for static and general nonlinear dynamic systems in continuous time \cite{KW:2000}, the theoretical advances of ES have gone much beyond such systems, encompassing discrete-time systems \cite{CKAL:2002}, stochastic systems \cite{MK:2009,LK:2010}, multivariable systems \cite{GKN:2012}, noncooperative games \cite{FKB:2012}, time delays \cite{OKT:2017,ZFO:2023} and even more general classes of infinite-dimensional systems governed by partial differential equations (PDEs) \cite{OK:2022_book_golden}. However, until recently there was no work that addressed ES with ETC---until \cite{VHPR:2023a,AUT:2025} which developed multi-variable ES algorithms based on perturbation-based (averaging-based) estimates of the model via ETC. Still, however, it remains an open problem to address concomitantly Nash equilibrium seeking (NES) for non-cooperative games \cite{DBB:1993,CHEN2023105510,MARTIROSYAN2024105883,NOORIGHANAVATIZADEH2023105642} using event-triggered versions \cite{KEK:2024,CAI2023105619,COUTINHO2024105678} of extremum seeking control \cite{K:2000_b,ZAGSK:2007,SK:2014,OK:2021,FOK:2023,SK:2023,ABDELGALIL2024105881}.}

In this paper and its companion conference version \cite{CDC:2024}, we advance such designs from the multi-input-single-output ES scenario to multi-input-multi-output NES scenario for games through the multi-variable ET-NES perspective. This is particularly challenging since in the literature for multi-input ES and ETC, most authors consider a centralized approach for the controller design (all vectors/matrices multiplying the control inputs need to be known), while games are clearly decentralized. We pursue a noncooperative game NES design where there are restrictions on sharing of information among the players in order to compute their tuning laws for obtaining the corresponding players' actions. Each player employing the proposed ET-NES distributed algorithm is improving its own payoff, irrespective of what the other players' actions are. We show that if all the players employ ET-NES algorithms they collectively converge to a Nash equilibrium. In other words, each of the players finds its optimal strategy, in an online fashion, even though they do not know the analytical forms of the payoff functions (neither the other players' nor their own) and neither have access to actions applied by the other players nor to payoffs achieved by the other players. While in \cite{CDC:2024} we have focused on duopoly games, here we address general $N$-player games, quadratic in the actions of the players. 
Our analysis employs time-scaling technique as well as proper sequencing steps of averaging for systems with discontinuous right-hand sides (due to the switching nature of the ETC actuation) and a Lyapunov function construction for the noncooperative result to prove closed-loop stability. A small neighborhood of the Nash equilibrium is achieved, while preserving performance even under limited bandwidth for the players' actions. We establish avoidance of the Zeno phenomenon
for both the average and the original system by proving the
existence of a minimal dwell-time between two successive event
triggers so that the tuning law does not have an infinite number of switchings or updates in any finite interval. A numerical example with an oligopoly game illustrates that the proposed ET-NES scheme differs from standard periodic sampled-data control \cite{KNTM:2013,HNW:2023,ZFO:2023}, since the event times (which result from the triggering condition and the system's state evolution) are in general only a (specific) subset of the sampling times and can be aperiodic. 


\section{$N$-Player Game with Quadratic Payoffs:\\ General Formulation}

As discussed earlier, game theory provides an important framework for mathematical modeling and analysis of scenarios involving different agents (players) where there is coupling in their actions, in the sense that their respective outcomes (outputs) $y_{i}(t) \in \mathbb{R}$ do not depend exclusively on their own actions/strategies (inputs) $\theta_{i}(t)\in \mathbb{R}$, with $i=1 ,\ldots, N$, but at least on a subset of others'. Moreover, defining $\theta := [\theta_1, \ldots, \theta_N]^{\top}$, each player's payoff function $J_{i}(\theta) : \mathbb{R}^{N} \to \mathbb{R}$ depends on the action $\theta_j$ of at least one other player $j$, $j\not= i$. An $N$-tuple of actions, $\theta^*$ is said to be in Nash equilibrium, if no player $i$ can improve his payoff  by unilaterally deviating from $\theta_i^*$, this being so for all $i$ \cite{Basar:1999}. 

Hence, we consider games where the payoff function of each player is quadratic, expressed as a strictly concave\footnote{By strict concavity, we mean $J_i(\theta)$ is strictly concave in $\theta_i$ for all $\theta_{-i}$, this being so for each $i=1,\ldots, N$, with $\theta_{-i}$ denoting the actions of the other players.} combination of their actions 
\begin{align}
J_{i}(\theta(t))=\frac{1}{2}\sum_{j=1}^{N}\sum_{k=1}^{N}H_{jk}^{i}\theta_{j}(t)\theta_{k}(t)+\sum_{j=1}^{N}h_{j}^{i}\theta_{j}(t)+c_{i}\,, \label{eq:Ji}
\end{align}  
where $\theta_{j}(t)  \in  \mathbb{R}$ is the decision variable of player $j$, $H_{jk}^{i}$, $h_{j}^{i}$, $c_{i}  \in  \mathbb{R}$ are constants, $H_{ii}^{i} < 0$, and 
\begin{align}
H_{jk}^{i}=H_{kj}^{i}\,, \quad \forall i,j,k \,. \label{eq:Hijk} 
\end{align}
Fig.~\ref{fig:blockDiagram_1} shows the closed-loop structure of the proposed ET-NES system to be designed.

Quadratic payoff functions are of particular interest in game theory, first because they constitute second-order approximations to other types of non-quadratic payoff functions, and second because they are analytically tractable, leading in general to closed-form equilibrium solutions which provide insight into the properties and features of the equilibrium solution concept under consideration \cite{Basar:1999}. 

For the sake of completeness, we provide here in mathematical terms, the definition of a Nash equilibrium $\theta^*=[\theta^*_1\,, \ldots\,,\theta_N^*]^{\top}$ in an $N$-player game:
\begin{align} \label{Nashcu}
J_i(\theta_i^*\,,\theta_{-i}^*) \geq J_i(\theta_i\,,\theta_{-i}^*)\,, \quad 
\forall \theta_i, \forall  i \in \{1\,, \ldots\,,N\}\,.
\end{align}
%
Hence, no player has an incentive to unilaterally deviate from its corresponding action from $\theta^*$.

In order to determine the Nash equilibrium solution in strictly concave quadratic games with $N$ players, where each action set is the entire real line, one should differentiate $J_{i}$ with respect to $\theta_{i}(t) \,, \forall i=1 ,\ldots, N$, setting the resulting expressions equal to zero, and solving the set of equations thus obtained. This set of equations, which also provides a sufficient condition due to strict concavity, is given by 
\begin{align}
\sum_{j =1}^{N}H_{ij}^{i}\theta_{j}^{*}+h_{i}^{i}=0\,,\quad i=1 ,\ldots, N\,, \label{eq:NE_v0}
\end{align} 
which can be written in the form of matrices as 
\begin{align}
\begin{bmatrix}
 H_{11}^{1} &  H_{12}^{1} & \hdots &  H_{1N}^{1} \\
 H_{21}^{2} &  H_{22}^{2} & \hdots &  H_{2N}^{2} \\
\vdots                       & \vdots                       &        & \vdots     \\
 H_{N1}^{N} &  H_{N2}^{N} & \hdots & H_{NN}^{N}   
\end{bmatrix}
\begin{bmatrix}
\theta_{1}^{*} \\
\theta_{2}^{*} \\
\vdots \\
\theta_{N}^{*} 
\end{bmatrix}
=-
\begin{bmatrix}
h_{1}^{1} \\
h_{2}^{2} \\
\vdots    \\
h_{N}^{N}   
\end{bmatrix}
\,,
\end{align}
which we re-write as 
\begin{align}
 H \theta^* = -h \label{eq:Htheta*h}
\end{align} 
from which we conclude that 
%
there exists only one Nash equilibrium at  $\theta^{*}=-H^{-1}h$, provided that $H$ is invertible. 
For more details, see \cite[Chapter 4]{Basar:1999}.

The \textit{objective} now is to design an extremum seeking-based algorithm to reach the Nash equilibrium using decentralized event-triggered policies.

\subsection{Continuous-time Nash Equilibrium Seeking}

Let us introduce for each $i$-th player the output signals 
\begin{align}
y_{i}(t)&=J_{i}(\theta(t))\,, \label{eq:yi_v0}
\end{align}
and introduce $\hat{\theta}_{i}(t)$ as an estimate of $\theta^{*}_{i}$ with \textit{estimation error} defined by
\begin{align}
\tilde{\theta}_{i}(t)&=\hat{\theta}_{i}(t)-\theta_{i}^{*}\,.\label{eq:tildeThetai}
\end{align}
Introduce additive-multiplicative dither signals $S_i(t)$ and $M_i(t)$, defined by 
\begin{align}
S_{i}(t)&=a_{i}\sin(\omega_{i}t) \,, \label{eq:Si} \\
M_{i}(t)&=\frac{2}{a_{i}}\sin(\omega_{i}t)\,, \label{eq:Mi}
\end{align}
with nonzero constant amplitudes $a_i>0$ at frequencies $\omega_i \neq \omega_j$. Here, the probing frequencies $\omega_{i}$'s can be selected as
\begin{align}
\omega_{i}=\omega_{i}'\omega \,, \quad i \in \left\{1,\ldots\,,N\right\}\,, \label{eq:omegai_event}
\end{align}
where $\omega$ is a positive constant and $\omega_{i}'$ is a rational number.

\begin{assum}\label{assumption3}
As considered in \cite{GKN:2012}, the probing frequencies are taken to satisfy
\begin{align}
\omega'_{i} 	\notin \left\{\omega'_{j}\,,~\frac{1}{2}(\omega'_{j}+\omega'_{k})\,,~\omega'_{j}+2\omega'_{k}\,,~\omega'_{k}\pm\omega'_{l}\right\}\,, \label{eq:omega_iNotIn}
\end{align}
for all $i$, $j$, $k$ and $l$.
\end{assum} 


\begin{figure*}[h!]
\begin{center}
\includegraphics[width=14.5cm]{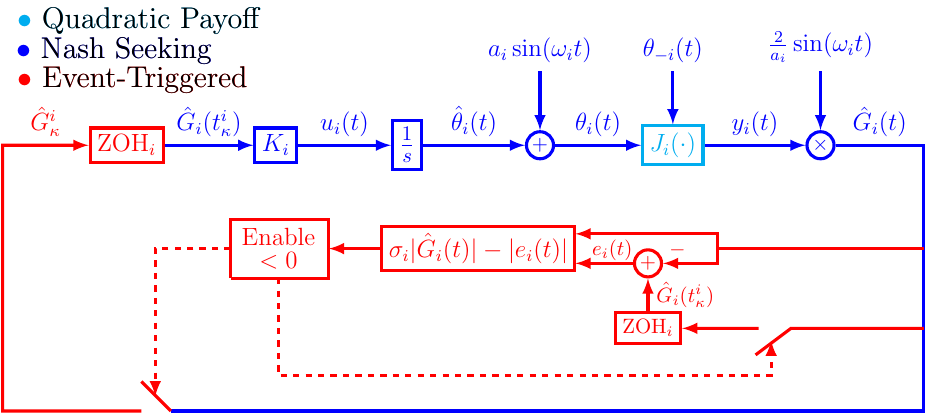}
\end{center}
\caption{Block diagram of the NES strategy through distributed event-triggered tuning policies of Definition~\ref{def:staticEvent} performed for each player.}
\label{fig:blockDiagram_1}
\end{figure*}

From Fig.~\ref{fig:blockDiagram_1}, the estimate $\hat{G}_{i}(t)$ of the unknown \textit{pseudo-gradient}\footnote{The \textit{pseudo-gradient} is essentially a generalized notion of gradient that provides directional information allowing the analysis of equilibrium points in complex strategic interactions of concave $N$-player games. Hence, given a vector-valued function $J(\theta): \mathbb{R}^{N} \to \mathbb{R}^{N}$, we define its \textit{pseudo-gradient} $G$ as the vector-valued function $G(\theta):=\left[\frac{\partial J_1}{\partial \theta_1}\,, \frac{\partial J_2}{\partial \theta_2}\,, \ldots\,, \frac{\partial J_N}{\partial \theta_N}\right]^T$. For more details, see \cite[Eq. (3.9)]{R:1965}.} of each payoff $J_i(\theta(t))$ is given by
\begin{align}
\hat{G}_{i}(t)&=M_{i}(t)y_{i}(t)\,. \label{eq:hatGi_nonaverage_1}
\end{align}
Plugging  (\ref{eq:Mi}) and (\ref{eq:yi_v0}) into (\ref{eq:hatGi_nonaverage_1}), we obtain 
\begin{align}
&\hat{G}_{i}(t)=\textcolor{blue}{\sum_{j=1}^{N}H_{ij}^{i}\tilde{\theta}_{j}(t)+\frac{2}{a_{i}}\sin(\omega_{i} t)\sum_{j=1}^{N}h_{j}^{i}\tilde{\theta}_{j}(t)} \nonumber \\
&\mathbb{+}\textcolor{blue}{\sum_{j=1}^{N} \sum_{k=1}^{N}\left[\frac{2\theta_{k}^{\ast}}{a_{i}}\sin(\omega_{i}t)\mathbb{-}\frac{a_{k}}{a_{i}}\cos\left[(\omega_{i}+\omega_{k})t\right]\right]H_{jk}^{i}\tilde{\theta}_{j}(t)}\nonumber \\
&\textcolor{blue}{\mathbb{+}\sum_{j=1}^{N} \sum_{\substack{k=1 \\ k\neq i}}^{N}\frac{a_{k}}{a_{i}}\cos\left[(\omega_{i}\mathbb{-}\omega_{k})t\right]H_{jk}^{i}\tilde{\theta}_{j}(t)} \nonumber \\
&\mathbb{+}\textcolor{red}{\sum_{j=1}^{N} \sum_{k=1}^{N}\frac{H_{ij}^{i}}{a_{i}}\sin(\omega_{i}t)\tilde{\theta}_{j}(t)\tilde{\theta}_{k}(t)}+\textcolor{magenta}{\sum_{j=1}^{N}H_{ij}^{i}\theta_{j}^{\ast}+h_{i}^{i}} \nonumber \\
&\mathbb{+}\frac{1}{2}\sum_{j=1}^{N} \sum_{k=1}^{N}H_{jk}^{i}\left\{-\frac{2a_{j}\theta_{k}^{\ast}}{a_{i}}\cos\left[(\omega_{i}\mathbb{+}\omega_{j})t\right] \right. \nonumber \\
&\mathbb{+}\frac{a_{j}a_{k}}{2a_{i}}\sin\left[(\omega_{i}\mathbb{+}\omega_{j}\mathbb{-}\omega_{k})t\right]\mathbb{+}\frac{a_{j}a_{k}}{2a_{i}}\sin\left[(\omega_{i}\mathbb{-}\omega_{j}\mathbb{+}\omega_{k})t\right]  \nonumber \\
&\mathbb{-}\frac{a_{j}a_{k}}{2a_{i}}\sin\left[(\omega_{i}\mathbb{+}\omega_{j}\mathbb{+}\omega_{k})t\right]\mathbb{-}\frac{a_{j}a_{k}}{2a_{i}}\sin\left[(\omega_{i}\mathbb{-}\omega_{j}\mathbb{-}\omega_{k})t\right] \nonumber \\
&\left.\mathbb{+}\frac{2\theta_{j}^{\ast}\theta_{k}^{\ast}}{a_{i}}\sin(\omega_{i}t)\right\}+\sum_{j=1}^{N} \sum_{\substack{k=1 \\ k\neq i}}^{N}\frac{a_{k}}{a_{i}}\cos\left[(\omega_{i}-\omega_{k})t\right]H_{jk}^{i}\theta_{j}^{\ast} \nonumber \\
&+\sum_{j=1}^{N}h_{j}^{i}\left\{-\frac{a_{j}}{a_{i}}\cos\left[(\omega_{i}+\omega_{j})t\right]+\frac{2\theta_{j}^{\ast}}{a_{i}}\sin(\omega_{i} t)\right\}\nonumber \\
&\mathbb{+}\sum_{\substack{j=1 \\ j\neq i}}^{N}\cos\left[(\omega_{i}-\omega_{j})t\right]\frac{a_{j}h_{j}^{i}}{a_{i}}+\frac{2c_{i}}{a_{i}}\sin(\omega_{i} t) \,. \label{eq:hatG_20250201_1}
\end{align}

From (\ref{eq:hatG_20250201_1}), it follows that the terms highlighted in \textcolor{blue}{blue} can be represented as a linear combination of the error variables $\tilde{\theta}_{j}(t)$, for all $j \in \{1\,, \ldots\,, N\}$, such that the time-varying coefficients corresponding to each $\tilde{\theta}_{j}(t)$ is given by
\begin{align}
&\mbox{\calligra H}_{~~ij}^{~~i}(t) \mathbb{=} H_{ij}^{i} \mathbb{+} \frac{2}{a_{i}}\sin(\omega_{i} t) h_{l}^{i} \mathbb{+} \sum_{k=1}^{N}\frac{2H_{lk}^{i}\theta_{k}^{\ast}}{a_{i}}\sin(\omega_{i}t) \nonumber \\
&\mathbb{-}\sum_{k=1}^{N}\frac{a_{k}H_{lk}^{i}}{a_{i}}\cos\left[(\omega_{i}\mathbb{+}\omega_{k})t\right]\mathbb{+} \sum_{k\neq i} \frac{a_{k}H_{lk}^{i}}{a_{i}}\cos\left[(\omega_{i}\mathbb{-}\omega_{k})t\right]. \label{eq:calligraHij}
\end{align}
The term in \textcolor{red}{red} is quadratic in $\tilde{\theta}(t)$ and, therefore, may be neglected in a local analysis \cite{AK:2003}. And, from (\ref{eq:NE_v0}), the term in \textcolor{magenta}{magenta} is equal to zero. Thus, by defining the time-varying disturbance
\begin{align}
&\Delta_{i}(t)=\frac{1}{2}\sum_{j=1}^{N} \sum_{k=1}^{N}H_{jk}^{i}\left\{-\frac{2a_{j}\theta_{k}^{\ast}}{a_{i}}\cos\left[(\omega_{i}\mathbb{+}\omega_{j})t\right] \right. \nonumber \\
&\mathbb{+}\frac{a_{j}a_{k}}{2a_{i}}\sin\left[(\omega_{i}\mathbb{+}\omega_{j}\mathbb{-}\omega_{k})t\right]\mathbb{+}\frac{a_{j}a_{k}}{2a_{i}}\sin\left[(\omega_{i}\mathbb{-}\omega_{j}\mathbb{+}\omega_{k})t\right]  \nonumber \\
&\mathbb{-}\frac{a_{j}a_{k}}{2a_{i}}\sin\left[(\omega_{i}\mathbb{+}\omega_{j}\mathbb{+}\omega_{k})t\right]\mathbb{-}\frac{a_{j}a_{k}}{2a_{i}}\sin\left[(\omega_{i}\mathbb{-}\omega_{j}\mathbb{-}\omega_{k})t\right] \nonumber \\
&\!\! \left.\mathbb{+}\frac{2\theta_{j}^{\ast}\theta_{k}^{\ast}}{a_{i}}\sin(\omega_{i}t)\right\}+\sum_{j=1}^{N} \sum_{\substack{k=1 \\ k\neq i}}^{N}\frac{a_{k}}{a_{i}}\cos\left[(\omega_{i}-\omega_{k})t\right]H_{jk}^{i}\theta_{j}^{\ast} \nonumber \\
&\mathbb{+}\sum_{j=1}^{N}h_{j}^{i}\left\{-\frac{a_{j}}{a_{i}}\cos\left[(\omega_{i}+\omega_{j})t\right]+\frac{2\theta_{j}^{\ast}}{a_{i}}\sin(\omega_{i} t)\right\}\nonumber \\
&\mathbb{+}\sum_{\substack{j=1 \\ j\neq i}}^{N}\cos\left[(\omega_{i}-\omega_{j})t\right]\frac{a_{j}h_{j}^{i}}{a_{i}}+\frac{2c_{i}}{a_{i}}\sin(\omega_{i} t)\,, \label{eq:Deltai}
\end{align}
the pseudo-gradient estimate (\ref{eq:hatG_20250201_1}) can be rewritten as
\begin{align}
&\hat{G}_{i}(t)=\sum_{j=1}^{N}\mbox{\calligra H}_{~~ij}^{~~i}(t)\tilde{\theta}_{j}(t)+\Delta_{i}(t) \,. \label{eq:hatG_20250206_1}
\end{align}

Defining the following time-varying matrix $\mbox{\calligra H}(t)$ and vector $\Delta(t)$,  
\begin{align}
\mbox{\calligra H}~~(t) &:= \begin{bmatrix}
													\mbox{\calligra H}_{~~11}^{~~1}(t) & \mbox{\calligra H}_{~~12}^{~~1}(t) & \ldots & \mbox{\calligra H}_{~~1N}^{~~1}(t) \\
													\mbox{\calligra H}_{~~21}^{~~2}(t) & \mbox{\calligra H}_{~~22}^{~~2}(t) & \ldots & \mbox{\calligra H}_{~~2N}^{~~2}(t) \\
													\vdots                         & \vdots                         &        & \vdots                         \\
													\mbox{\calligra H}_{~~N1}^{~~N}(t) & \mbox{\calligra H}_{~~N2}^{~~N}(t) & \ldots & \mbox{\calligra H}_{~~NN}^{~~N}(t) \\
												 \end{bmatrix} \,, \label{eq:calligraH} \\
						 \Delta(t) &:= \begin{bmatrix}
													\Delta_{1}(t) \,, 
													\Delta_{2}(t) \,,
													\ldots\,,
													\Delta_{N}(t)
												 \end{bmatrix}^{\top} \,, \label{eq:Delta}
\end{align}
and $\hat{G}(t)=[\hat{G}_{1}(t),,\hat{G}_{2}(t),,\ldots,,\hat{G}_{N}(t)]^{\top}$ $\in \mathbb{R}^{N}$ as the vector of the pseudo-gradient estimate, we can express (\ref{eq:hatG_20250206_1}) in the next compact form 
\begin{align}
\hat{G}(t)&=\mbox{\calligra H}~~(t)\tilde{\theta}(t)+\Delta(t)\,, \label{eq:hatG_20240302_2}
\end{align}
where $\mbox{\calligra H }(t)$, determined by (\ref{eq:calligraH}), is a time-varying matrix with average value equal to $H$, while $\Delta(t)$, given by (\ref{eq:Delta}), is a time-varying vector of zero mean.

On the other hand, from the time-derivative of (\ref{eq:tildeThetai}) and the NES scheme depicted in Fig.~\ref{fig:blockDiagram_1}, the dynamics that governs $\hat{\theta}(t)$, as well as $\tilde{\theta}(t)$, are given by
\begin{align}
\frac{d\tilde{\theta}(t)}{dt}&=\frac{d\hat{\theta}(t)}{dt}=u(t) \label{eq:dtildeThetadt_20250206_1}\,, 
\end{align}
with $u(t)= [u_{1}(t)\,,u_{2}(t)\,,\ldots\,,u_{N}(t)]^{\top}$ where each $u_{j}(t)$, for all $j \in \{1\,, \ldots \,, N\}$, is a decentralized ES law to be designed. Moreover, by using (\ref{eq:dtildeThetadt_20250206_1}), the time-derivative of (\ref{eq:hatG_20240302_2}) is given by 
\begin{align}
\frac{d\hat{G}(t)}{dt}&=\mbox{\calligra H}~~(t)u(t)+\frac{d\mbox{\calligra H}~~(t)}{dt}\tilde{\theta}(t)+\frac{d \Delta(t)}{dt}\,, \label{eq:dhatGdt_20250206_1}
\end{align}
where, from (\ref{eq:calligraH}) and (\ref{eq:Delta}), it is easy to verify that $\dfrac{d\mbox{\calligra H}~~(t)}{dt}$ as well as $\dfrac{d \Delta(t)}{dt}$ are of zero mean. 
 
For all $t\geq 0$, the  continuous-time proportional feedback law
\begin{align}
u_{i}(t)=K_{i}\hat{G}_{i}(t) \,, \quad \forall t\geq 0 \label{eq:U_continuous}
\end{align}
is a stabilizing tuning law with the  gain
\begin{align}
K=\text{diag}\left\{K_{1}\,,K_{2}\,,\ldots\,,K_{N}\right\}\,,
\end{align}
being chosen such that $KH$ is Hurwitz. This can be easily obtained for any $K>0$ and $H$ in (\ref{eq:Htheta*h}) satisfying the next assumption.
\begin{assum} \label{assum:SDD} The \textit{unknown }matrix $H$ is strictly diagonally dominant \cite{FKB:2012}, {\it i.e.},
\begin{align}
\sum_{j \neq i} |H_{ij}^{i}| < |H_{ii}^{i}|, \quad \forall i \in \{1, \dots, N\}. 
\end{align}
\end{assum}
Under Assumption~\ref{assum:SDD}, $H$ is invertible and the Nash equilibrium $\theta^{\ast}$ exists and is unique.

Here, updates of the pseudo-gradient estimate are triggered only for a given  sequence of time instants $(t_{\kappa})_{\kappa\in\mathbb{N}}$ determined by an event-generator. This generator is designed to maintain stability and robustness. The task is orchestrated by a monitoring mechanism that triggers updates when the difference between the current output value and its previously computed value at time $t_\kappa$ exceeds a predefined threshold, determined by a constructed triggering condition \cite{HJT:2012}. It is important to note that in conventional sampled-data implementations, execution times are evenly spaced in time, with $t_{\kappa+1}=t_{\kappa}+ h$, where $h>0$ is a known constant, for all $\kappa\in \mathbb{N}$. However, in an event-triggered scheme, sampling times may occur aperiodically, as desired. 

\subsection{Emulation of the Continuous-Time\\ Nash Equilibrium  Seeking Design}

Let the $i$-th tuning law be
\begin{align}
u_i(t)=K_{i}\hat{G}_{i}(t^{i}_{\kappa}) \,, \quad \forall t\in[t^{i}_{\kappa},t^{i}_{\kappa+1})\,, \quad \kappa\in\mathbb{N}\,, \label{eq:U_event}
\end{align}
we can define the error function $e_{i} : \mathbb{R}  \mapsto\mathbb{R} $ as the difference between the $i$-th player current state variable and its last broadcasted value as
\begin{align}
e_{i}(t):=\hat{G}_{i}(t^{i}_{\kappa})-\hat{G}_{i}(t) \,, \quad \forall t \in \lbrack t^{i}_{\kappa}\,, t^{i}_{\kappa+1}) \,, \quad \kappa\in \mathbb{N} \,. \label{eq:e_event}
\end{align}

Now, using  (\ref{eq:e_event}) and  (\ref{eq:U_event}), for all $t\in[t^{i}_{\kappa},t^{i}_{\kappa+1})$,  the $i$-th event-triggered tuning law is rewritten as  
\begin{align}
u_i(t)=K_{i}\hat{G}_{i}(t)+K_{i}e_{i}(t) \,, \quad \forall t \in \lbrack t^{i}_{\kappa}\,, t^{i}_{\kappa+1}) \,, \quad \kappa\in \mathbb{N} \,.  \label{eq:U_event24}
\end{align}

Consequently, defining $e(t)= [e_{1}(t)\,,e_{2}(t)\,,\ldots\,,e_{N}(t)]^{\top}$ and using the pseudo-gradient estimate (\ref{eq:hatG_20240302_2}), the dynamics (\ref{eq:dtildeThetadt_20250206_1}) and (\ref{eq:dhatGdt_20250206_1}), and the decentralized event-triggered tuning law (\ref{eq:U_event24}), the closed-loop system governing $\hat{G}(t)$ and $\tilde{\theta}(t)$ is given by
\begin{align}
\frac{d\hat{G}(t)}{dt}&=\mbox{\calligra H}~~(t)K\hat{G}(t)+\mbox{\calligra H}~~(t)Ke(t) \nonumber \\
&\quad+\frac{d\mbox{\calligra H}~~(t)}{dt}\tilde{\theta}(t)+\frac{d \Delta(t)}{dt}\,, \label{eq:dhatGdt_20250206_2} \\
\frac{d\tilde{\theta}(t)}{dt}&=K\hat{G}(t)+Ke(t) \nonumber \\
&= K\mbox{\calligra H}~~(t)\tilde{\theta}(t)+K\Delta(t)+Ke(t)\,, \label{eq:dtildeThetadt_20250206_2}
\end{align}
for all $t\in[t_{\kappa},t_{\kappa+1})$, where $t_{\kappa}=\min\{t^{1}_{\kappa}\,,\ldots\,,t^{N}_{\kappa}\}$ and $t_{\kappa+1}=\min\{t^{1}_{\kappa+1}\,,\ldots\,,t^{N}_{\kappa+1}\}$ (for more details, see \cite{TC:2014}).

The closed-loop system described by (\ref{eq:dhatGdt_20250206_2}) and (\ref{eq:dtildeThetadt_20250206_2}) highlights a crucial point: while the product $\mbox{\calligra H}~~(t)K$ on averaging sense forms a Hurwitz matrix, the convergence to the equilibrium $\hat{G}\equiv0$ and $\tilde{\theta}\equiv0$ is not guaranteed due to the presence of the error vector $e(t)$ and the time-varying term $\Delta(t)$ and their derivatives. However, the system does exhibit Input-to-State Stability (ISS) concerning the error vector $e(t)$ and such time-varying disturbances. Additionally, it is important to note that the disturbances $\Delta(t)$ and $\frac{d\Delta(t)}{dt}$ as well as the time-varying matrix $\frac{d\mbox{\calligra H}~~(t)}{dt}$ possess zero mean values.

In the next two sections, we introduce a distributed static event-triggering mechanism for NES, as outlined in Definitions \ref{def:staticEvent} and \ref{def:averageStaticEvent}. In this mechanism, players abstain from sharing their state information with others, thus creating a noncooperative game environment. 
Within this setup of data transmission, players autonomously determine their broadcast actions based solely on their individual information. This
mechanism represents a fusion of distributed event-triggered
data transmission with an extremum-seeking tuning system \cite{VHPR:2023a,AUT:2025}.

\section{Distributed Event-Triggered Tuning for Nash Equilibrium Seeking} 

The scheme involves each state variable being monitored by a separate player. In this scenario, no single sensor has access to the complete state vector $\hat{G}(t)=[\hat{G}_{1}(t), G_{2}(t) \,,\ldots\,, G_{N}(t)]^{\top}$ required by the event-triggered extremum seeking mechanisms in \cite{VHPR:2022,VHPR:2023b,VHPR:2023a,AUT:2025}. Here, the aim is to come up with a scheme whereby individual players can assess independently the difference between the current value of their state and its last broadcasted value to determine when to execute a new triggering update of the estimate of the pseudo-gradient. 


Definition~\ref{def:staticEvent} below articulates the utilization of small parameters $\sigma_{i} \in (0,1)$, along with the errors $e_{i}$ representing the disparity between a player's current state variable and its last broadcasted value, and measurements of the pseudo-gradient estimate $\hat{G}_{i}$. These components are employed to construct the ``NES Static-Triggering Condition''. This approach involves updating the estimate of the pseudo-gradient, at triggering times, by means of the ZOH actuators in order to obtain the tuning laws (\ref{eq:U_event}), as depicted in Fig.~\ref{fig:blockDiagram_1}. This ensures the asymptotic stability of the closed-loop system.

\begin{definition}[\small{NES Static-Triggering Condition}] \label{def:staticEvent}
The ET-NES with static-triggering condition consists of two components:
\begin{enumerate}
	\item The set of increasing sequences of time $I=\{I_{1}\,,\ldots\,,I_{N}\}$ such that $I_{i}=\{t^{i}_{0}\,, t^{i}_{1}\,, t^{i}_{2}\,,\ldots\}$, with $t^{i}_{0}=0$, for all $i\in \{1,\ldots,N\}$, generated under the following rules:
		\begin{itemize}
			\item If $\left\{t \in\mathbb{R}^{+}\mathbb{:}~ t\mathbb{>}t^{i}_{\kappa} ~ \mathbb{\wedge} ~ \sigma_{i}|\hat{G}_{i}(t)| -|e_{i}(t)| < 0 \right\} = \emptyset$, then the set of the times of the events is $I_{i}=\{t^{i}_{0}\,, t^{i}_{1}\,, \ldots, t^{i}_{\kappa}\}$.
			\item If $\left\{t \in\mathbb{R}^{+}\mathbb{:}~ t\mathbb{>}t^{i}_{\kappa} ~ \mathbb{\wedge} ~ \sigma_{i}|\hat{G}_{i}(t)| -|e_{i}(t)| < 0 \right\} \neq \emptyset$, the next event time is given by
				\begin{align}
					\!\!\!\! t^{i}_{\kappa+1}&\mathbb{=}\inf\left\{t \in\mathbb{R}^{+}\mathbb{:}~ t\mathbb{>}t^{i}_{\kappa} ~ \mathbb{\wedge} ~ \sigma_{i}|\hat{G}_{i}(t)| -|e_{i}(t)| < 0 \right\}, \label{eq:tk+1_event}
				\end{align}
				$\!\!\!$consisting of the static event-triggering mechanism.
		\end{itemize}
	\item The $i$-th player tuning law using the pseudo-gradient updates at the triggering instants is (\ref{eq:U_event}).
\end{enumerate}  
\end{definition}

\section{Closed-Loop System for Time-Scaled\\ Triggering Mechanism} \label{sec:cls}

\subsection{Rescaling of Time} \label{sec:timeScale}

Now, we  introduce a suitable time scale to carry out  the stability analysis of the closed-loop system. From (\ref{eq:omegai_event}), we can note that the dither frequencies  (\ref{eq:Si}) and (\ref{eq:Mi}), as well as their combinations, are both rational. Furthermore, there exists a time period $T$ such that 
\begin{align}
T&= 2\pi \times \text{LCM}\left\{\frac{1}{\omega_{i}}\right\}\,, \quad \forall i \left\{1\,,2\,,\ldots\,,N\right\}\,, \label{eq:T}
\end{align}
 with LCM denoting the least common multiple.  Hence, it is possible to define the time-scale for the dynamics  (\ref{eq:dtildeThetadt_20250206_1}) and (\ref{eq:dhatGdt_20250206_1}) with the transformation $\bar{t}=\omega t$, where 
\begin{align}
\omega&:=\frac{2\pi}{T}\,. \label{eq:omega_event_1}
\end{align}
 Consequently, the system of equations (\ref{eq:dtildeThetadt_20250206_1}) and (\ref{eq:dhatGdt_20250206_1}) can be rewritten, $\forall t \in \lbrack t_{\kappa}\,, t_{\kappa+1})$ and $\kappa\in \mathbb{N}$, as 
\begin{align}
\frac{d\hat{G}(\bar{t})}{d\bar{t}}&=\frac{1}{\omega}\mbox{\calligra H}~~(\bar{t})K\hat{G}(\bar{t})+\frac{1}{\omega}\mbox{\calligra H}~~(\bar{t})Ke(\bar{t}) \nonumber \\
&\quad+\frac{d\mbox{\calligra H}~~(\bar{t})}{d\bar{t}}\tilde{\theta}(\bar{t})+\frac{1}{\omega}\frac{d \Delta(\bar{t})}{d\bar{t}}\,, \label{eq:dhatGdt_20250206_3} \\
\frac{d\tilde{\theta}(\bar{t})}{d\bar{t}}&=\frac{1}{\omega}K\hat{G}(\bar{t})+\frac{1}{\omega}Ke(\bar{t}) \nonumber \\
&= \frac{1}{\omega}K\mbox{\calligra H}~~(\bar{t})\tilde{\theta}(\bar{t})+\frac{1}{\omega}K\Delta(\bar{t})+\frac{1}{\omega}Ke(\bar{t})\,.\label{eq:dtildeThetadt_20250206_3}
\end{align}
At this point, we can implement a suitable averaging mechanism within the transformed time scale $\bar{t}$ based on the dynamics (\ref{eq:dhatGdt_20250206_3}) and (\ref{eq:dtildeThetadt_20250206_3}). Despite the non-periodicity of the triggering events and discontinuity on the right-hand sides, the closed-loop system maintains its periodicity over time due to the periodic probing and demodulation signals. This unique characteristic allows for the application of the averaging results established by Plotnikov \cite{P:1979} to this particular setup---see Appendix~\ref{appendix_plotnikov}. 


\subsection{Average Closed-Loop System}

Defining the augmented state as follows
\begin{align}
X(\bar{t}):=\begin{bmatrix} \hat{G}(\bar{t})\,, \tilde{\theta}(\bar{t})\end{bmatrix}^{\top}\,,
\end{align}
the system (\ref{eq:dhatGdt_20250206_3})--(\ref{eq:dtildeThetadt_20250206_3}) is reduced to
\begin{align}
\dfrac{dX(\bar{t})}{d\bar{t}}&=\dfrac{1}{\omega}\mathcal{F}\left(\bar{t},X,\dfrac{1}{\omega}\right)\,. \label{eq:dotX_event}
\end{align}
The system (\ref{eq:dotX_event}) features a small parameter $1/\omega$ as well as a $T$-periodic function $\mathcal{F}\left(\bar{t},X,\dfrac{1}{\omega}\right)$ in $\bar{t}$. Therefore, the averaging theorem \cite{P:1979} can be applied to  $\mathcal{F}\left(\bar{t},X,\dfrac{1}{\omega}\right)$ at $\displaystyle \lim_{\omega\to \infty}\dfrac{1}{\omega}=0$. The averaging method allows for determining in what sense the behavior of a constructed average autonomous system approximates the behavior of the non-autonomous system (\ref{eq:dotX_event}). Of course, it can be inferred intuitively that in instances where the response of a system is significantly slower than its excitation, the response will predominantly be dictated by the average characteristics of the excitation.

By employing the averaging computation to (\ref{eq:dotX_event}), we derive the following average system
\begin{align}
\dfrac{dX_{\text{av}}(\bar{t})}{d\bar{t}}&=\dfrac{1}{\omega}\mathcal{F}_{\text{av}}\left(X_{\text{av}}\right) \,, \label{eq:dotXav_event_1} \\
\mathcal{F}_{\text{av}}\left(X_{\text{av}}\right)&=\dfrac{1}{T}\int_{0}^{T}\mathcal{F}\left(\xi,X_{\text{av}},0\right)d\xi
\,.  \label{eq:mathcalFav_event}
\end{align}
Therefore, ``freezing'' the average states of $\hat{G}(\bar{t})$, $e(\bar{t})$, and $\tilde{\theta}(\bar{t})$ in (\ref{eq:dhatGdt_20250206_3})--(\ref{eq:dtildeThetadt_20250206_3}), the averaging terms are
\begin{align}
\mbox{\calligra H}^{~~\rm{av}}(\bar{t})&= \frac{1}{T}\int_{0}^{T}\mbox{\calligra H}~~(\xi)d\xi=H \,, \label{eq:calligraHav} \\
\frac{d\mbox{\calligra H}^{~~\rm{av}}(\bar{t})}{d \bar{t}}&= \frac{1}{T}\int_{0}^{T}\frac{d\mbox{\calligra H}~~(\xi)}{d \bar{t}}d\xi=0 \,, \label{eq:calligradHavdt} \\
\Delta^{\rm{av}}(\bar{t})&= \frac{1}{T}\int_{0}^{T}\Delta(\xi)d\xi=0 \,, \label{eq:Deltaav} \\
\frac{d\Delta^{\rm{av}}(\bar{t})}{d \bar{t}}&= \frac{1}{T}\int_{0}^{T}\frac{d\Delta(\xi)}{d \bar{t}}d\xi=0 \,, \label{eq:dDeltaavdt}
\end{align}
and we get, for all $\bar t\in [\bar t_k, \bar t_{k+1})$: 
\begin{align}
\frac{d\hat{G}_{\text{av}}(\bar{t})}{d\bar{t}}&=\frac{1}{\omega}KH\hat{G}_{\text{av}}(\bar{t})+\frac{1}{\omega}KHe_{\text{av}}(\bar{t})\,, \label{eq:dotHatGav_event_1} \\
\frac{d\tilde{\theta}_{\text{av}}(\bar{t})}{d\bar{t}}&=\frac{1}{\omega}KH\tilde{\theta}_{\text{av}}(\bar{t})+\frac{1}{\omega}Ke_{\text{av}}(\bar{t})\,. \label{eq:dotTildeThetaAv_event_1} \\
e_{\text{av}}(\bar{t})&=\hat{G}_{\text{av}}(\bar{t}_{k})-\hat{G}_{\text{av}}(\bar{t})\,, \label{eq:Eav_event_1} \\
\hat{G}_{\text{av}}(\bar{t})&= H\tilde{\theta}_{\text{av}}(\bar{t})\,,\label{eq:hatGav_event_1}
\end{align}
recalling that the matrix $HK$ is Hurwitz. Thus, it is evident from (\ref{eq:dotHatGav_event_1}) that the ISS relationship of $\hat{G}_{\rm{av}}(\bar{t})$ with respect to the average measurement error $e_{\rm{av}}(\bar{t})$ holds. Thus, we can introduce the following ``Average Static-Triggering Condition'' for the average system.

\begin{definition}[\small{Average Static-Triggering Condition}] \label{def:averageStaticEvent} The average event-triggered condition consists of two components:
\begin{enumerate}
	\item The set of increasing sequences of time $I=\{I_{1}\,,\ldots\,,I_{N}\}$ such that $I_{i}=\{\bar{t}^{i}_{0}\,, \bar{t}^{i}_{1}\,, \bar{t}^{i}_{2}\,,\ldots\}$, with $\bar{t}^{i}_{0}=0$, for all $i\in \{1,\ldots,N\}$, generated under the following rules:
		\begin{itemize}
			\item If $\left\{\bar{t} \in\mathbb{R}^{+}\mathbb{:}~ \bar{t}\mathbb{>}\bar{t}^{i}_{\kappa} ~ \mathbb{\wedge} ~ \sigma_{i}|\hat{G}_{i}^{\rm{av}}(\bar{t})| \mathbb{-}|e_{i}^{\rm{av}}(\bar{t})| \mathbb{<}  0 \right\} = \emptyset$, then the set of the times of the events is $I_{i}=\{\bar{t}^{i}_{0}\,, \bar{t}^{i}_{1}\,, \ldots, \bar{t}^{i}_{\kappa}\}$.
			\item If $\left\{\bar{t} \in\mathbb{R}^{+}\mathbb{:}~ \bar{t}\mathbb{>}\bar{t}^{i}_{\kappa} ~ \mathbb{\wedge} ~ \sigma_{i}|\hat{G}_{i}^{\rm{av}}(\bar{t})| \mathbb{-}|e_{i}^{\rm{av}}(\bar{t})| \mathbb{<} 0 \right\} \neq \emptyset$, the next event time is given by
				\begin{align}
					\!\!\!\!\!\!\!\! \bar{t}^{i}_{\kappa+1}&=\mathbb{\inf}\left\{\bar{t} \in\mathbb{R}^{+}\mathbb{:}~ \bar{t}\mathbb{>}\bar{t}^{i}_{\kappa} ~ \mathbb{\wedge} ~ \sigma_{i}|\hat{G}_{i}^{\rm{av}}(\bar{t})| \mathbb{-}|e_{i}^{\rm{av}}(\bar{t})| \mathbb{<} 0 \right\}\!, \label{eq:tk+1_event_av}
				\end{align}
				$\!\!\!$consisting of the static event-triggering mechanism.
		\end{itemize}
		\item The $i$-th player tuning law using the average pseudo-gradient updates at the triggering instants is 
		\begin{align}
			u_{i}^{\rm av}(\bar{t})=K_{i}\hat{G}^{\rm av}_{i}(\bar{t}^{i}_{k}) \,, \label{eq:U_MD4}
		\end{align}
		for all $\bar{t} \in \lbrack \bar{t}_{k}\,, \bar{t}_{k+1}\phantom{(}\!\!)$, $k\in \mathbb{N}$.
\end{enumerate} 
\end{definition}

We claim that the corresponding event-triggering mechanisms introduced in Definitions \ref{def:staticEvent} and \ref{def:averageStaticEvent} guarantee the asymptotic stabilization of $\hat{G}_{\text{av}}(\bar{t})$ and, consequently, that of $\tilde{\theta}_{\text{av}}(\bar{t})$. Since $H$ is invertible, both $\hat{G}_{\text{av}}(\bar{t})$ and $\tilde{\theta}_{\text{av}}(\bar{t})$ converge to the origin according to the averaging theory \cite{K:2002}.


\section{Stability Analysis}\label{ETESNC_unknownH*}

The next theorem guarantees the local asymptotic stability of the ET-NES system employing static event-triggered execution, as depicted in Fig.~\ref{fig:blockDiagram_2}$.1$---see Appendix~\ref{appendix_plotnikov}.


\begin{theorem} \label{thm:NETESC_2}
Consider the closed-loop average dynamics of the pseudo-gradient estimate (\ref{eq:dotHatGav_event_1}), the average error vector (\ref{eq:Eav_event_1}), the average static event-triggering mechanism in Definition \ref{def:averageStaticEvent}, and Assumptions \ref{assumption3} and \ref{assum:SDD}. For $\omega>0$ sufficiently large, defined in (\ref{eq:omega_event_1}), the equilibrium $\hat{G}_{\text{av}}(t)=0$ is locally exponentially stable and $\tilde{\theta}_{\text{av}}(t)$ converges exponentially to zero. In particular, for
the non-average system (\ref{eq:dtildeThetadt_20250206_2}), depicted in Fig.~\ref{fig:blockDiagram_1}, there exist constants   $m\,,M_{\theta}>0$ such that
\begin{align}
\|\theta(t)-\theta^{\ast}\|&\leq M_{\theta}\exp(-mt)\|\theta(0) - \theta^{\ast}\|+\mathcal{O}\left(a+\frac{1}{\omega}\right), \label{eq:normTheta_thm2} 
\end{align}
where $a=\sqrt{\sum_{i=1}^{n}a_{i}^{2}}$, with $a_i$ defined in \eqref{eq:Si} and the constants $m$ and $M_{\theta}$ depending on the triggering parameters $\sigma_i$. In addition, there exists a lower bound  $\tau^{\ast}$ for the inter-execution interval $t_{k+1}-t_{k}$, for all $k \in \mathbb{N}$,  preventing the Zeno behavior.
\end{theorem} 
\textit{Proof:} The proof of the theorem is divided into two parts: stability analysis and avoidance of Zeno behavior.

\begin{flushleft}
\underline{\it A. Stability Analysis}
\end{flushleft}

Consider the following candidate Lyapunov function for the average system \eqref{eq:dotHatGav_event_1}: 
\begin{align}
V_{\text{av}}(\bar{t})=\hat{G}^{\top}_{\text{av}}(\bar{t})P\hat{G}_{\text{av}}(\bar{t}) \,,\, P=P^{\top} >0. \label{eq:Vav_event_pf2}
\end{align}  
Since $HK$ is Hurwitz, given $Q=Q^T>0$ there exists $P=P^T$ such that the Lyapunov equation is $K^{T}H^{T}P+PHK=-Q$ and, therefore, the time derivative of (\ref{eq:Vav_event_pf2}) is given by
\begin{align}
\frac{dV_{\text{av}}(\bar{t})}{d\bar{t}}&=-\frac{1}{\omega}\hat{G}_{\text{av}}^{\top}(\bar{t})Q\hat{G}_{\text{av}}(\bar{t})+\frac{1}{\omega}e_{\text{av}}^{\top}(\bar{t})H^{\top}K^{\top}P\hat{G}_{\text{av}}(\bar{t}) \nonumber \\
&\quad+\frac{1}{\omega}\hat{G}_{\text{av}}^{\top}(\bar{t})PKHe_{\text{av}}(\bar{t})\,, \label{eq:dotVav_event_1_pf2}
\end{align}
whose  upper bound satisfies
\begin{small}
\begin{align}
\frac{dV_{\text{av}}(\bar{t})}{d\bar{t}}&\!\leq\!\mathbb{-}\frac{\lambda_{\min}(Q)}{\omega}\|\hat{G}_{\text{av}}(\bar{t})\|^{2}\mathbb{+}
\frac{ 2\|PKH\|}{\omega}\|e_{\text{av}}(\bar{t})\| \|\hat{G}_{\text{av}}(\bar{t})\|. \label{eq:dotVav_event_2_pf2}
\end{align}
\end{small}

In the proposed event-triggering mechanism, the average update law is (\ref{eq:tk+1_event_av})--(\ref{eq:U_MD4}), and  $u^{\text{av}}_{i}(\bar{t})$ is held constant between two consecutive events. Therefore, the norm of the average measurement error $e^{\rm{av}}(\bar{t})$ is upper bounded by
\begin{align}
\|e_{\text{av}}(\bar{t})\|& = \sqrt{\sum_{j=1}^{N}|e_{i}^{\text{av}}(\bar{t})|^{2}} \leq \sqrt{\sum_{j=1}^{N}\sigma_{i}|\hat{G}_{i}^{\text{av}}(\bar{t})|^{2}}	\nonumber \\
													& \leq \bar{\sigma}\sqrt{\sum_{j=1}^{N}|\hat{G}_{i}^{\text{av}}(\bar{t})|^{2}}= \bar{\sigma}\|\hat{G}_{\text{av}}(\bar{t})\|\,, \label{eq:eAv_upperBound}
\end{align}
with 
\begin{align}
\bar{\sigma} = \max_{i\in \{1,\ldots , N\}}\left\{\sigma_{1},\sigma_{2},\ldots, \sigma_{N}\right\}\,.
\end{align}
Now, plugging (\ref{eq:eAv_upperBound}) into (\ref{eq:dotVav_event_2_pf2}), we obtain 
\begin{align}
\frac{dV_{\text{av}}(\bar{t})}{d\bar{t}}&\leq-\frac{\lambda_{\min}(Q)}{\omega}\left(1-\frac{ 2\|PKH\|\bar{\sigma}}{\lambda_{\min}(Q)}\right)\|\hat{G}_{\text{av}}(\bar{t})\|^{2}\,. \label{eq:dotVav_event_4_pf2}
\end{align}

By using the Rayleigh-Ritz Inequality \cite{K:2002}, we get 
\begin{align}
\lambda_{\min}(P)\|\hat{G}_{\text{av}}(\bar{t})\|^{2}\leq V_{\text{av}}(\bar{t}) \leq \lambda_{\max}(P)\|\hat{G}_{\text{av}}(\bar{t})\|^{2}\,, \label{eq:Rayleigh-Ritz_pf2}
\end{align}
 and the following upper bound for (\ref{eq:dotVav_event_4_pf2}): 
\begin{align}
\frac{dV_{\text{av}}(\bar{t})}{d\bar{t}}&\leq-\frac{1}{\omega}\frac{\lambda_{\min}(Q)}{\lambda_{\max}(P)}\left(1-\frac{ 2\|PKH\|\bar{\sigma}}{\lambda_{\min}(Q)}\right)V_{\text{av}}(\bar{t}) \,, \label{eq:dotVav_event_5_pf2}
\end{align}
with $V_{\text{av}}(\bar{t})>0$ and $\dfrac{dV_{\text{av}}(\bar{t})}{d\bar{t}}<0$, for all $\bar{\sigma}<\dfrac{\lambda_{\min}(Q)}{ 2\|PHK\|}$. For instance, if we choose $\bar{\sigma}=\frac{\lambda_{\min}(Q)}{ 2\|PHK\|}\hat{\sigma}$, where $\hat{\sigma} \in (0,1)$, and defining 
\begin{align}
\bar{\sigma}&=\frac{\lambda_{\min}(Q)}{ 2\|PHK\|}\hat{\sigma}\,, \quad  \hat{\sigma} \in (0,1)\,, \quad \alpha = \frac{\lambda_{\min}(Q)}{\lambda_{\max}(P)}\,, \label{eq:hatSigma_20240307_1}
\end{align}
inequality (\ref{eq:dotVav_event_5_pf2}) simply becomes 
\begin{align}
\frac{dV_{\text{av}}(\bar{t})}{d\bar{t}}&\leq-\frac{\alpha\left(1-\hat{\sigma}\right)}{\omega}V_{\text{av}}(\bar{t}) \,. \label{eq:dotVav_20240307_1}
\end{align}

Then, using \cite[Comparison Lemma]{K:2002} an upper bound $\bar{V}_{\text{av}}(\bar{t})$ for $V_{\text{av}}(\bar{t})$ according to 
\begin{align}
V_{\text{av}}(\bar{t})\leq \bar{V}_{\text{av}}(\bar{t}) \,, \quad \forall \bar{t}\in \lbrack \bar{t}_{k},\bar{t}_{k+1}\phantom{(}\!\!) \,, \label{eq:VavBarVav_1_pf2}
\end{align}
 is given by the solution of the equation
\begin{align}
\frac{d\bar{V}_{\text{av}}(\bar{t})}{d\bar{t}}=-\frac{\alpha\left(1-\hat{\sigma}\right)}{\omega}\bar{V}_{\text{av}}(\bar{t})\,, \quad \bar{V}_{\text{av}}(\bar{t}_{k})=V_{\text{av}}(\bar{t}_{k})\,.
\end{align}
In other words, $ \forall \bar{t}\in \lbrack \bar{t}_{k},\bar{t}_{k+1}\phantom{(}\!\!)$,
\begin{align}
\bar{V}_{\text{av}}(\bar{t})=\exp\left(-\frac{\alpha\left(1-\hat{\sigma}\right)}{\omega}\bar{t}\right)V_{\text{av}}(\bar{t}_{k})\,, \label{eq:_pf2}
\end{align}
and inequality (\ref{eq:VavBarVav_1_pf2}) is rewritten as
\begin{align}
V_{\text{av}}(\bar{t})\leq \exp\left(-\frac{\alpha\left(1-\hat{\sigma}\right)}{\omega}\bar{t}\right)V_{\text{av}}(\bar{t}_{k}) \,, \quad \forall \bar{t}\in \lbrack \bar{t}_{k},\bar{t}_{k+1}\phantom{(}\!\!)
\,. \label{eq:VavBarVav_2_pf2}
\end{align}

By defining $\bar{t}_{k}^{+}$ and $\bar{t}_{k}^{-}$ as the right and left limits of $\bar{t}=\bar{t}_{k}$, respectively, it easy to verify that $V_{\text{av}}(\bar{t}_{k+1}^{-})\leq \exp\left(-\dfrac{\alpha\left(1-\hat{\sigma}\right)}{\omega}(\bar{t}_{k+1}^{-}-\bar{t}_{k}^{+})\right)V_{\text{av}}(\bar{t}_{k}^{+})$. Since $V_{\text{av}}(\bar{t})$ is continuous, $V_{\text{av}}(\bar{t}_{k+1}^{-})=V_{\text{av}}(\bar{t}_{k+1})$ and $V_{\text{av}}(\bar{t}_{k}^{+})=V_{\text{av}}(\bar{t}_{k})$, and therefore,
\begin{align}
    V_{\text{av}}(\bar{t}_{k+1})\leq \exp\left(-\frac{\alpha\left(1-\hat{\sigma}\right)}{\omega}(\bar{t}_{k+1}-\bar{t}_{k})\right)V_{\text{av}}(\bar{t}_{k})\,. \label{METES_eq:mmd_1_s}
\end{align}
Hence, for any $\bar{t}\geq 0$ in $ \bar{t}\in \lbrack \bar{t}_{k},\bar{t}_{k+1}\phantom{(}\!\!)$, $k \in \mathbb{N}$, one has 
\begin{align}
    V_{\text{av}}(\bar{t})&\leq \exp\left(-\frac{\alpha\left(1-\hat{\sigma}\right)}{\omega}(\bar{t}-\bar{t}_{k})\right) V_{\text{av}}(\bar{t}_{k}) \nonumber \\
    &\leq \exp\left(-\frac{\alpha\left(1-\hat{\sigma}\right)}{\omega}(\bar{t}-\bar{t}_{k})\right)  \nonumber \\
		&\quad \times\exp\left(-\frac{\alpha\left(1-\hat{\sigma}\right)}{\omega}(\bar{t}_{k}-\bar{t}_{k-1})\right)V_{\text{av}}(\bar{t}_{k-1}) \nonumber \\
    &\leq \ldots \leq \nonumber \\
    &\leq \exp\left(-\frac{\alpha\left(1-\hat{\sigma}\right)}{\omega}(\bar{t}\mathbb{-}\bar{t}_{k})\right) \nonumber \\
		&\quad \times \prod_{i=1}^{i=k}\exp\left(-\frac{\alpha\left(1-\hat{\sigma}\right)}{\omega}(\bar{t}_{i}-\bar{t}_{i-1})\right)V_{\text{av}}(\bar{t}_{i-1}) \nonumber \\
    &=\exp\left(-\frac{\alpha\left(1-\hat{\sigma}\right)}{\omega}\bar{t}\right) V_{\text{av}}(0)\,, \quad \forall \bar{t}\geq 0\,. \label{METES_eq:VavBarVav_2_pf2}
\end{align}
Now, by lower bounding the left-hand side and upper bounding the right-hand side of (\ref{METES_eq:VavBarVav_2_pf2}) with their counterparts in Rayleigh-Ritz inequality (\ref{eq:Rayleigh-Ritz_pf2}), we obtain: 
\begin{align}
\lambda_{\mathbb{\min}}(P)\|\hat{G}_{\text{av}}(\bar{t})\|^{2}&\mathbb{\leq} \mathbb{\exp}\left(\!\!\mathbb{-}\frac{\alpha\left(1\mathbb{-}\hat{\sigma}\right)}{\omega}\bar{t}\right)\!\!\lambda_{\max}(P)\|\hat{G}_{\text{av}}(0)\|^{2}\!\!. \label{eq:VavBarVav_3_pf2}
\end{align}
Then,
\begin{align}
\|\hat{G}_{\text{av}}(\bar{t})\|^{2}&\mathbb{\leq} \mathbb{\exp}\left(\!\!\mathbb{-}\frac{\alpha\left(1\mathbb{-}\hat{\sigma}\right)}{\omega}\bar{t}\right)\frac{\lambda_{\max}(P)}{\lambda_{\min}(P)}\|\hat{G}_{\text{av}}(0)\|^{2} \nonumber \\
&\mathbb{=}\left[\mathbb{\exp}\left(\!\!\mathbb{-}\frac{\alpha\left(1\mathbb{-}\hat{\sigma}\right)}{2\omega}\bar{t}\right)\sqrt{\frac{\lambda_{\max}(P)}{\lambda_{\min}(P)}}\|\hat{G}_{\text{av}}(0)\|\right]^{2}\!\!, \label{eq:VavBarVav_4_pf2}
\end{align}
and
\begin{align}
\|\hat{G}_{\text{av}}(\bar{t})\|\leq\exp\left(-\frac{\alpha\left(1-\hat{\sigma}\right)}{2\omega}\bar{t}\right)\sqrt{\frac{\lambda_{\max}(P)}{\lambda_{\min}(P)}}\|\hat{G}_{\text{av}}(0)\|\,. \label{eq:normHatGav_1_pf2}
\end{align}
Since $H$ is invertible, from (\ref{eq:hatGav_event_1}), $\tilde{\theta}_{\rm{av}}(\bar{t})\mathbb{=}H^{-1}\hat{G}_{\rm{av}}(\bar{t})$ and $\|\tilde{\theta}_{\rm{av}}(\bar{t})\|\mathbb{=}\|H^{-1}\hat{G}_{\rm{av}}(\bar{t})\|\mathbb{\leq} \|H^{-1}\| \|\hat{G}_{\rm{av}}(\bar{t})\|$. Consequently,
\begin{align}
\|H^{-1}\|^{-1}\|\tilde{\theta}_{\rm{av}}(\bar{t})\|\leq  \|\hat{G}_{\rm{av}}(\bar{t})\|\,. \label{eq:normG_l}
\end{align}
From (\ref{eq:hatGav_event_1}), we also conclude 
\begin{align}
 \|\hat{G}_{\rm{av}}(\bar{t})\| \leq \|H\|\|\tilde{\theta}_{\rm{av}}(\bar{t})\|\,. \label{eq:normG_h}
\end{align}
Thus, by using (\ref{eq:normG_l}) and (\ref{eq:normG_h}), inequality (\ref{eq:normHatGav_1_pf2}) can be rewritten with respect to $\tilde{\theta}_{\rm{av}}(\bar{t})$ as
\begin{align}
\|\tilde{\theta}_{\text{av}}(\bar{t})\|&\leq\exp\left(-\frac{\alpha\left(1-\hat{\sigma}\right)}{2\omega}\bar{t}\right)\sqrt{\frac{\lambda_{\max}(P)}{\lambda_{\min}(P)}} \nonumber \\
&\quad \times \|H^{-1}\|\|H\|\|\tilde{\theta}_{\text{av}}(0)\| \,. \label{eq:normTildeThetaAv_1_pf2}
\end{align}
Since (\ref{eq:dotTildeThetaAv_event_1}) has a discontinuous right-hand side, but it is also $T$-periodic in $t$, and noting that the average system with state $\tilde{\theta}_{\text{av}}(\bar{t})$ is asymptotically stable according to (\ref{eq:normTildeThetaAv_1_pf2}), we can invoke the averaging theorem in \cite[Theorem~2]{P:1979} to conclude that
\begin{align}
\|\tilde{\theta}(t)-\tilde{\theta}_{\text{av}}(t)\|\leq\mathcal{O}\left(\frac{1}{\omega}\right)\,.
\end{align}
By applying the triangle inequality \cite{A:1957}, we also obtain:
\begin{align}
\|\tilde{\theta}(t)\|&\leq\|\tilde{\theta}_{\text{av}}(t)\|+\mathcal{O}\left(\frac{1}{\omega}\right)\nonumber \\
&\leq \exp\left(-\frac{\alpha\left(1-\hat{\sigma}\right)}{2}t\right)M_{\theta}\|\tilde{\theta}_{\text{av}}(0)\|+\mathcal{O}\left(\frac{1}{\omega}\right)\!. \label{eq:nomrTildeTheta_pf2_gold}
\end{align}
Now, from (\ref{eq:Si}) and Fig.~\ref{fig:blockDiagram_1}, we can verify  that
\begin{align}
\theta(t)-\theta^{\ast}=\tilde{\theta}(t)+S(t)\,, \label{eq:theta_3_event_pf2}
\end{align}
where $S(t):=[S_{1}(t)\,,\ldots\,,S_{N}(t)]^{\top}$ and whose the Euclidean norm satisfies
\begin{align}
&\|\theta(t)-\theta^{\ast}\|=\|\tilde{\theta}(t)+S(t)\| \leq \|\tilde{\theta}(t)\|+\|S(t)\| \nonumber \\
&\leq \exp\left(-\frac{\alpha\left(1-\hat{\sigma}\right)}{2}t\right)M_{\theta}\|\theta(0) - \theta^{\ast}\|+\mathcal{O}\left(a+\frac{1}{\omega}\right)\,, \label{eq:theta_4_event_pf2}
\end{align}
leading to inequality (\ref{eq:normTheta_thm2}), for $m=\frac{\alpha(1-\hat{\sigma})}{2}$ and $M_{\theta}=\sqrt{\frac{\lambda_{\max}(P)}{\lambda_{\min}(P)}} \|H^{-1}\|\|H\|$.

\begin{flushleft}
\textcolor{black}{\underline{\it B. Avoidance of Zeno Behavior}}
\end{flushleft}

Since the average closed-loop system consists of (\ref{eq:dotHatGav_event_1}), with the event-triggering mechanism  (\ref{eq:tk+1_event_av}) and the average tuning law (\ref{eq:U_MD4}), we can conclude that $\|e_{\text{av}}(\bar{t})\| \leq \bar{\sigma}\|\hat{G}_{\text{av}}(\bar{t})\|$, resulting in \begin{align}
\bar{\sigma}\|\hat{G}_{\text{av}}(\bar{t})\|^{2}-\|e_{\text{av}}(\bar{t})\|\|\hat{G}_{\text{av}}(\bar{t})\|\geq 0\,. \label{ineq:interEvents_1_static}
\end{align}
By using the Peter-Paul inequality \cite{W:1971}, $cd\leq \frac{c^2}{2\epsilon}+\frac{\epsilon d^2}{2}$ for all $c,d,\epsilon>0$, with $c=\|e_{\rm{av}}(\bar{t})\|$, $d=\|\hat{G}_{\rm{av}}(\bar{t})\|$ and $\epsilon=\bar{\sigma}$, the inequality (\ref{ineq:interEvents_1_static}) is lower bounded by
\begin{align}
&\bar{\sigma} \|\hat{G}_{\text{av}}(\bar{t})\|^{2}-\|e_{\text{av}}(\bar{t})\|\|\hat{G}_{\text{av}}(\bar{t})\|\geq \nonumber \\
&\bar{\sigma} \|\hat{G}_{\text{av}}(\bar{t})\|^{2}-\left(\frac{\bar{\sigma}}{2}\|\hat{G}_{\rm{av}}(\bar{t})\|^2+\frac{1}{2\bar{\sigma}}\|e_{\rm{av}}(\bar{t})\|^2\right)\nonumber \\
&= q\|\hat{G}_{\rm{av}}(\bar{t})\|^{2}-p\|e_{\rm{av}}(\bar{t})\|^2\,,\label{ineq:interEvents_2_static_pf2}
\end{align}
where $q=\frac{\bar{\sigma}}{2}$ and $p=\frac{1}{2\bar{\sigma}}$. 
%
In \cite{G:2014}, it is shown that a lower bound for the inter-execution interval is given by the time duration it takes for the function
\begin{align}
\phi_{\rm{av}}(\bar{t})=\sqrt{\frac{p}{q}}\frac{\|e_{\rm{av}}(\bar{t})\|}{\|\hat{G}_{\rm{av}}(\bar{t})\|} \label{eq:phi_1_static_pf2}
\end{align}
to go from 0 to 1. The time-derivative of (\ref{eq:phi_1_static_pf2}) is 
\begin{align}
\frac{d\phi_{\rm{av}}(\bar{t})}{d\bar{t}}&=\sqrt{\frac{p}{q}}\frac{1}{\|e_{\rm{av}}(\bar{t})\|\|\hat{G}_{\rm{av}}(\bar{t})\|}\left[e_{\rm{av}}^{\top}(\bar{t})\frac{de_{\rm{av}}(\bar{t})}{d\bar{t}} \right. \nonumber \\
&\quad\left.-\hat{G}_{\rm{av}}^{\top}(\bar{t})\frac{d\hat{G}_{\rm{av}}(\bar{t})}{d\bar{t}}\left(\frac{\|e_{\rm{av}}(\bar{t})\|}{\|\hat{G}_{\rm{av}}(\bar{t})\|}\right)^2\right]\,. \label{eq:dotPhi_1_static_pf2}
\end{align}
Now, plugging (\ref{eq:dotHatGav_event_1}) and (\ref{eq:Eav_event_1}) into (\ref{eq:dotPhi_1_static_pf2}), the following expression is arrived at:  
\begin{align}
&\frac{d\phi_{\rm{av}}(\bar{t})}{d\bar{t}}\mathbb{=}\frac{1}{\omega}\sqrt{\frac{p}{q}}\frac{1}{\|e_{\rm{av}}(\bar{t})\|\|\hat{G}_{\rm{av}}(\bar{t})\|} \nonumber \\
& \mathbb{\times}\left\{\mathbb{-}e_{\rm{av}}^{\top}(\bar{t})KHe_{\rm{av}}(\bar{t}) \mathbb{-}e_{\rm{av}}^{\top}(\bar{t})KH\hat{G}_{\rm{av}}(\bar{t})+\right.\nonumber \\
&\left. \mathbb{-}\left[\hat{G}_{\rm{av}}^{\top}(\bar{t})KH\hat{G}_{\rm{av}}(\bar{t})\mathbb{+}\hat{G}_{\rm{av}}^{\top}(\bar{t})KHe_{\rm{av}}(\bar{t})\right]\!\!\left(\frac{\|e_{\rm{av}}(\bar{t})\|}{\|\hat{G}_{\rm{av}}(\bar{t})\|}\right)^{\!\!\! 2}\!\right\}\!\!. \label{eq:dotPhi_1_1_static_pf2}
\end{align}
Then, the following estimate holds:
\begin{align}
\frac{d\phi_{\rm{av}}(\bar{t})}{d\bar{t}}&\leq
\frac{\|KH\|}{\omega}\sqrt{\frac{p}{q}}\left(1+\frac{\|e_{\rm{av}}(\bar{t})\|}{\|\hat{G}_{\rm{av}}(\bar{t})\|}\right)^2\,. \label{eq:dotPhi_2_static_pf2}
\end{align}

Hence,  using (\ref{eq:phi_1_static_pf2}), inequality (\ref{eq:dotPhi_2_static_pf2}) is rewritten as
\begin{align}
\frac{d\phi_{\text{av}}(\bar{t})}{d\bar{t}}&\leq\frac{\|KH\|}{\omega}\sqrt{\frac{q}{p}}\left(\sqrt{\frac{p}{q}}+\phi_{\text{av}}(\bar{t})\right)^{2}\,. \label{eq:dotPhi_3_static_pf2}
\end{align}
By invoking \cite[Comparison Lemma]{K:2002}, the following upper bound  for $\phi_{\text{av}}(\bar{t})$ is obtained 
\begin{align}
\phi_{\text{av}}(\bar{t})\leq \tilde{\phi}_{\text{av}}(\bar{t}) \,, \quad \phi_{\text{av}}(0)= \tilde{\phi}_{\text{av}}(0)=0 \,, \label{eq:tildePhi_v2}
\end{align}
where $\tilde{\phi}_{\text{av}}(\bar{t})$ is the solution of the equation
\begin{align}
\frac{d\tilde{\phi}\textcolor{black}{_{\text{av}}}(\bar{t})}{d\bar{t}}&=\frac{\|KH\|}{\omega}\sqrt{\frac{q}{p}}\left(\sqrt{\frac{p}{q}}+\tilde{\phi}_{\text{av}}(\bar{t})\right)^{2}\,. \label{eq:dotTildePhi_v2}
\end{align}
The solution of (\ref{eq:dotTildePhi_v2}), with the initial condition $\tilde{\phi}_{\text{av}}(0) = 0$, is simply 
\begin{align}
\tilde{\phi}_{\text{av}}(\bar{t}) = \dfrac{\sqrt{\dfrac{p}{q}}}{1 - \dfrac{\|KH\|}{\omega} \dfrac{q}{p} \bar{t}} - \sqrt{\frac{p}{q}}.  \label{eq:phi-bart_v2}
\end{align}
Since $\phi_{\text{av}}(t)$ in (\ref{eq:phi_1_static_pf2}) is an average version of $\phi(t)=\sqrt{\frac{p}{q}}\frac{|e(\bar{t})|}{|\hat{G}(\bar{t})|}$, by invoking \cite[Theorem 2]{P:1979}, one can write down 
\begin{align}
|\phi(t)-\tilde{\phi}_{\text{av}}(t)|\leq\mathcal{O}\left(\frac{1}{\omega}\right)\,.
\end{align}
By using the triangle inequality \cite{A:1957}, one has
\begin{align}
\phi(t)&\leq\ \phi_{\text{av}}(t)+\mathcal{O}\left(\frac{1}{\omega}\right) \leq \tilde{\phi}_{\text{av}}(t)+\mathcal{O}\left(\frac{1}{\omega}\right)\! \nonumber \\
&=\dfrac{\sqrt{\frac{p}{q}}}{1 - \dfrac{\|KH\|}{\omega} \dfrac{q}{p} t} - \sqrt{\dfrac{p}{q}} +\mathcal{O}\left(\frac{1}{\omega}\right). \label{eq:phi_t}
\end{align}
Now, defining 
\begin{align}
\hat{\phi}(t):=\dfrac{\sqrt{\dfrac{p}{q}}}{1 - \|KH\| \dfrac{q}{p} t} - \sqrt{\dfrac{p}{q}} +\mathcal{O}\left(\frac{1}{\omega}\right)\,, \label{eq:hatPhi_t}
\end{align}
a lower bound on the inter-execution interval for the original system (or its average version) is given by the time it takes for the function (\ref{eq:hatPhi_t}), or (\ref{eq:phi-bart_v2}), to go from 0 to 1.  Considering (\ref{eq:hatPhi_t}), this lower bound is at least equal to 
\begin{align} \label{cataflan_novalgina}
\tau^{\ast}&=\dfrac{1}{\|KH\|}\frac{1}{\bar{\sigma}^2}\frac{1-\mathcal{O}(1/\omega)}{1+1/\bar{\sigma}-\mathcal{O}(1/\omega)}\,.
\end{align} 
Therefore, the Zeno behavior is avoided not only for the average closed-loop system in time $\bar{t}$ but also for the original system in time $t$ since $\bar{t}={\omega}t$ is only a time compression (dilation) for $\omega$ sufficiently large (small), which means that $\bar{\tau}^*={\omega}\tau$ will still be a finite number, establishing a minimum switching time to rule out any Zeno behavior in time $t\in\mathbb{R}_{+}$ or $\bar{t}\in\mathbb{R}_{+}$. 
\hfill $\square$

\section{Simulation Results}

This section presents simulation results in order to illustrate the distributed NES-based event-triggering scheme. The investigated system captures a noncooperative game involving four firms operating in an oligopoly market structure. 

These firms engage in competition to maximize their profits   
%
%
\begin{align}
J_{1}(t)&=\frac{1}{2}\theta^{\top}(t)H_{1}\theta(t)+h_{1}^{\top}\theta(t)+c_{1}\,, \label{eq:J1_20240311} \\
J_{2}(t)&=\frac{1}{2}\theta^{\top}(t)H_{2}\theta(t)+h_{2}^{\top}\theta(t)+c_{2}\,, \label{eq:J2_20240311} \\
J_{2}(t)&=\frac{1}{2}\theta^{\top}(t)H_{3}\theta(t)+h_{3}^{\top}\theta(t)+c_{3}\,, \label{eq:J3_20240311} \\
J_{4}(t)&=\frac{1}{2}\theta^{\top}(t)H_{4}\theta(t)+h_{4}^{\top}\theta(t)+c_{4}\,, \label{eq:J4_20240311}
\end{align}
by setting the prices $\theta_{i}(t)$ of their respective products without sharing any information among the players (see \cite[section 4.6]{Basar:1999} and \cite{FKB:2012} for background). We denote the vector of prices by $\theta(t):=[\theta_{1}(t) \ \theta_{2}(t) \ \theta_{3}(t) \ \theta_{4}(t)]^{\top}\in \mathbb{R}^{4}$ and define the remaining parameters as follows: 
\begin{small}
\begin{align}
&H_{1}= \frac{1}{(R_{2}R_{3}R_{4}+R_{1}R_{3}R_{4}+R_{1}R_{2}R_{4}+R_{1}R_{2}R_{3})}\times \nonumber \\
&\times\begin{bmatrix} -2(R_{3}R_{4}+R_{2}R_{4}+R_{2}R_{3}) &R_{3}R_{4} & R_{2}R_{4} & R_{2}R_{3} \\ 
R_{3}R_{4} & 0 & 0 & 0 \\
R_{2}R_{4} & 0 & 0 & 0 \\
R_{2}R_{3} & 0 & 0 & 0 \end{bmatrix} \,, \\ 
&h_{1} = \frac{1}{(R_{2}R_{3}R_{4}+R_{1}R_{3}R_{4}+R_{1}R_{2}R_{4}+R_{1}R_{2}R_{3})} \times \nonumber \\
&\times\begin{bmatrix} m_{1}(R_{3}R_{4}+R_{2}R_{4}+R_{2}R_{3})+S_{d}R_{2}R_{3}R_{4}  \\ 
-m_{1}R_{3}R_{4}  \\
-m_{1}R_{2}R_{4}  \\
-m_{1}R_{2}R_{3}  \end{bmatrix}\,, \\
&c_{1} = -\frac{m_{1}S_{d}R_{2}R_{3}R_{4}}{(R_{2}R_{3}R_{4}+R_{1}R_{3}R_{4}+R_{1}R_{2}R_{4}+R_{1}R_{2}R_{3})}\,,
\end{align}
\end{small}
\begin{small}
\begin{align}
&H_{2}= \frac{1}{(R_{2}R_{3}R_{4}+R_{1}R_{3}R_{4}+R_{1}R_{2}R_{4}+R_{1}R_{2}R_{3})}\times \nonumber \\
&\times\begin{bmatrix} 0 & R_{3}R_{4} & 0 & 0 \\ 
 R_{3}R_{4} & -2(R_{3}R_{4}+R_{1}R_{4}+R_{1}R_{3}) & R_{1}R_{4} & R_{1}R_{3} \\
0 & R_{1}R_{4} & 0 & 0 \\
0 & R_{1}R_{3} & 0 & 0 \end{bmatrix} \,, \\ 
&h_{2} = \frac{1}{(R_{2}R_{3}R_{4}+R_{1}R_{3}R_{4}+R_{1}R_{2}R_{4}+R_{1}R_{2}R_{3})} \times \nonumber \\
&\times\begin{bmatrix} -m_{2}R_{3}R_{4} \\ 
m_{2}(R_{3}R_{4}+R_{1}R_{4}+R_{1}R_{3})+S_{d}R_{1}R_{3}R_{4}  \\
-m_{2}R_{1}R_{4}  \\
-m_{2}R_{1}R_{3}  \end{bmatrix}\,, \\
&c_{2} = -\frac{m_{2}S_{d}R_{1}R_{3}R_{4}}{(R_{2}R_{3}R_{4}+R_{1}R_{3}R_{4}+R_{1}R_{2}R_{4}+R_{1}R_{2}R_{3})}\,,
\end{align}
\end{small}
\begin{small}
\begin{align}
&H_{3}= \frac{1}{(R_{2}R_{3}R_{4}+R_{1}R_{3}R_{4}+R_{1}R_{2}R_{4}+R_{1}R_{2}R_{3})}\times \nonumber \\
&\times \begin{bmatrix} 0 & 0 & R_{2}R_{4} & 0 \\ 
 0 & 0 &  R_{1}R_{4} & 0 \\
R_{2}R_{4} & R_{1}R_{4} & -2(R_{2}R_{4}+R_{1}R_{4}+R_{1}R_{2}) & R_{1}R_{2} \\
0 & 0 & R_{1}R_{2} & 0 \end{bmatrix} \,, \\ 
&h_{3} = \frac{1}{(R_{2}R_{3}R_{4}+R_{1}R_{3}R_{4}+R_{1}R_{2}R_{4}+R_{1}R_{2}R_{3})}\times \nonumber \\
&\times\begin{bmatrix} -m_{3}R_{2}R_{4} \\ 
-m_{3}R_{1}R_{4}  \\
m_{3}(R_{2}R_{4}+R_{1}R_{4}+R_{1}R_{2})+S_{d}R_{1}R_{2}R_{4}  \\
-m_{3}R_{1}R_{2}  \end{bmatrix}\,, \\
&c_{3} = -\frac{m_{3}S_{d}R_{1}R_{2}R_{4}}{(R_{2}R_{3}R_{4}+R_{1}R_{3}R_{4}+R_{1}R_{2}R_{4}+R_{1}R_{2}R_{3})}\,,
\end{align}
\end{small}
\begin{small}
\begin{align}
&H_{4}= \frac{1}{(R_{2}R_{3}R_{4}+R_{1}R_{3}R_{4}+R_{1}R_{2}R_{4}+R_{1}R_{2}R_{3})}\times \nonumber \\
&\times \begin{bmatrix} 0 & 0 & 0 & R_{2}R_{3} \\ 
 0 & 0 &  0 & R_{1}R_{3} \\
0 & 0 & 0 & R_{1}R_{2} \\
R_{2}R_{3} & R_{1}R_{3} & R_{1}R_{2} & -2(R_{2}R_{3}+R_{1}R_{3}+R_{1}R_{2})\end{bmatrix} \,, \\ 
&h_{4} = \frac{1}{(R_{2}R_{3}R_{4}+R_{1}R_{3}R_{4}+R_{1}R_{2}R_{4}+R_{1}R_{2}R_{3})}\times \nonumber \\
&\times\begin{bmatrix} -m_{4}R_{2}R_{3} \\ 
-m_{4}R_{1}R_{3}  \\
-m_{4}R_{1}R_{2} \\
m_{4}(R_{2}R_{3}+R_{1}R_{3}+R_{1}R_{2})+S_{d}R_{1}R_{2}R_{3}  \end{bmatrix}\,, \\
&c_{4} = -\frac{m_{4}S_{d}R_{1}R_{2}R_{3}}{(R_{2}R_{3}R_{4}+R_{1}R_{3}R_{4}+R_{1}R_{2}R_{4}+R_{1}R_{2}R_{3})}\,,
\end{align}
\end{small}

\medskip

where $m_{1}$, $m_{2}$, $m_{3}$ and $m_{4}$ are marginal costs, $S_{d}$ is the total consumer demand, and $R_{1}$, $R_{2}$, $R_{3}$ and $R_{4}$ represent the resistance that consumers have toward buying a given product. This resistance may be due to quality or brand image considerations---the most desirable products have lowest resistance. The payoff functions (\ref{eq:J1_20240311})--(\ref{eq:J4_20240311}) are in the form (\ref{eq:Ji}). Therefore, the Nash equilibrium $\theta^{\ast} = [\theta_{1}^{\ast}\,, \theta_{2}^{\ast}\,, \theta_{3}^{\ast}\,, \theta_{4}^{\ast}]^{\top}$ satisfies (\ref{eq:Htheta*h}) with $H$ and $h$ given by (\ref{eq:H_20240311}) and (\ref{eq:h_20240311}), respectively.
\begin{small}
\begin{figure*}
\begin{align}
&H= \frac{1}{(R_{2}R_{3}R_{4}+R_{1}R_{3}R_{4}+R_{1}R_{2}R_{4}+R_{1}R_{2}R_{3})}\times \nonumber \\
&\times\begin{bmatrix} -2(R_{3}R_{4}+R_{2}R_{4}+R_{2}R_{3}) &R_{3}R_{4} & R_{2}R_{4} & R_{2}R_{3} \\ 
 R_{3}R_{4} & -2(R_{3}R_{4}+R_{1}R_{4}+R_{1}R_{3}) & R_{1}R_{4} & R_{1}R_{3} \\
R_{2}R_{4} & R_{1}R_{4} & -2(R_{2}R_{4}+R_{1}R_{4}+R_{1}R_{2}) & R_{1}R_{2} \\
R_{2}R_{3} & R_{1}R_{3} & R_{1}R_{2} & -2(R_{2}R_{3}+R_{1}R_{3}+R_{1}R_{2})\end{bmatrix} \,, \label{eq:H_20240311} \\
&h=\frac{1}{(R_{2}R_{3}R_{4}+R_{1}R_{3}R_{4}+R_{1}R_{2}R_{4}+R_{1}R_{2}R_{3})}\begin{bmatrix} m_{1}(R_{3}R_{4}+R_{2}R_{4}+R_{2}R_{3})+S_{d}R_{2}R_{3}R_{4}  \\ 
m_{2}(R_{3}R_{4}+R_{1}R_{4}+R_{1}R_{3})+S_{d}R_{1}R_{3}R_{4}  \\
m_{3}(R_{2}R_{4}+R_{1}R_{4}+R_{1}R_{2})+S_{d}R_{1}R_{2}R_{4}  \\
m_{4}(R_{2}R_{3}+R_{1}R_{3}+R_{1}R_{2})+S_{d}R_{1}R_{2}R_{3}  \end{bmatrix}\,. \label{eq:h_20240311}
\end{align}
\end{figure*}
\end{small}
In addition, according to Assumption~\ref{assum:SDD}, the matrix $H$ in (\ref{eq:H_20240311}) is strictly diagonally dominant. Hence, Nash equilibrium $\theta^{\ast}$ exists and is unique since strictly diagonally dominant matrices are nonsingular by Levy-Desplanques' Theorem \cite{T:1949,HJ:1985}. Moreover, $H^{i}_{ij}=H_{j}^{ji}$, and hence $H$ is a negative definite symmetric matrix and the nooncoperative game belongs to the class of games known as \textit{potential games} \cite{MS:1996} or games that are strategically equivalent to strictly concave team problems \cite{Basar:1999}. \newpage

We have simulated the decentralized event-triggered implementation of static noncooperative game with four firms in an oligopoly market structure following the design procedure presented in this paper. The parameters of the plant are the same as those in \cite{FKB:2012}: with initial conditions  $\hat{\theta}_{1}(0)=52$, $\hat{\theta}_{2}(0)=40.93$, $\hat{\theta}_{3}(0)=33.5$, $\hat{\theta}_{4}(0)=35.09$; $S_{d} = 100$, $R_{1}=0.15$, $R_{2}=0.30$, $R_{3}=0.60$, $R_{4}=1$, $m_{1}=30$, $m_{2}=30$, $m_{3}=25$, $m_{4}=20$, which, according to (\ref{eq:Htheta*h}), yield the unique Nash equilibrium 
\begin{align} 
\theta^{*}&=\begin{bmatrix}42.8818 & 40.9300 & 37.8363 & 35.0874\end{bmatrix}^{\top}\,, \label{eq:theta_NE} \\
J^{*}&=\begin{bmatrix}524.0208 & 293.4217 & 238.4846 & 209.6584\end{bmatrix}^{\top}\,. \label{eq:J_NE}
\end{align}
Moreover, the parameters employed for the tuning policies were:  $a_{1}=a_{2}=a_{3}=a_{4}=0.05$, $K_{1}=6$, $K_{2}=18$, $K_{3}=10$, $K_{4}=24$, $\omega_{1}=30$, $\omega_{2}=24$, $\omega_{3}=44$, $\omega_{4}=36$,   $\sigma_{1}=0.65$, $\sigma_{2}=0.55$, $\sigma_{1}=0.75$ and $\sigma_{4}=0.45$.

\begin{figure*}[h!]
	\centering
	\subfigure[Aperiodic behavior of the players' tuning laws $u_i(t)=K_{i}\hat{G}_{i}(t^{i}_{\kappa})$, based on a piecewise-constant/ZOH use of the estimate of the pseudo-gradient $\hat{G}_{i}(t)$. \label{fig:u}]{\includegraphics[width=6.27cm]{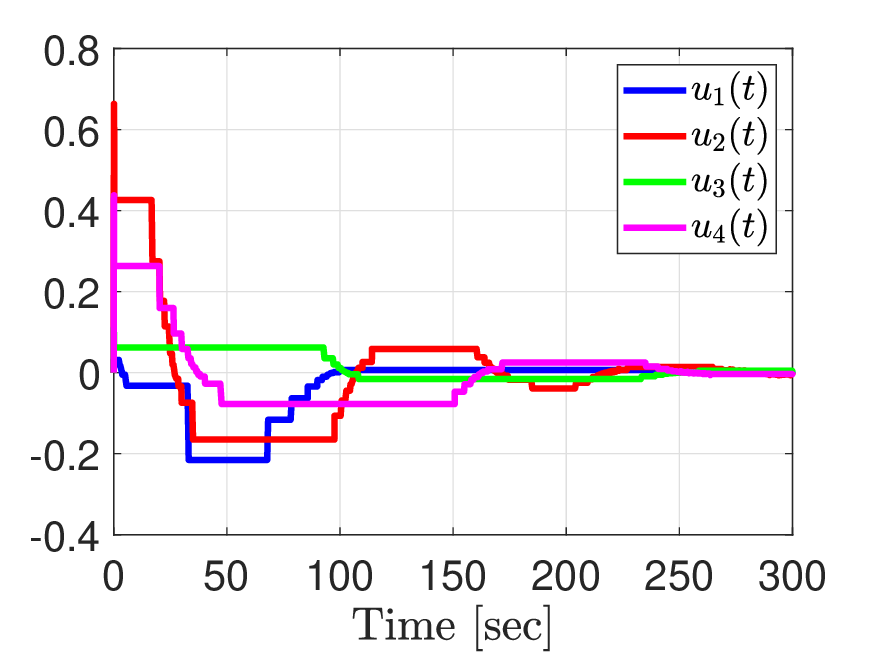}}
	~
	\subfigure[Time evolution of the triggering updates of the estimates of the pseudo-gradients. \label{fig:controlUpdate}]{\includegraphics[width=6.27cm]{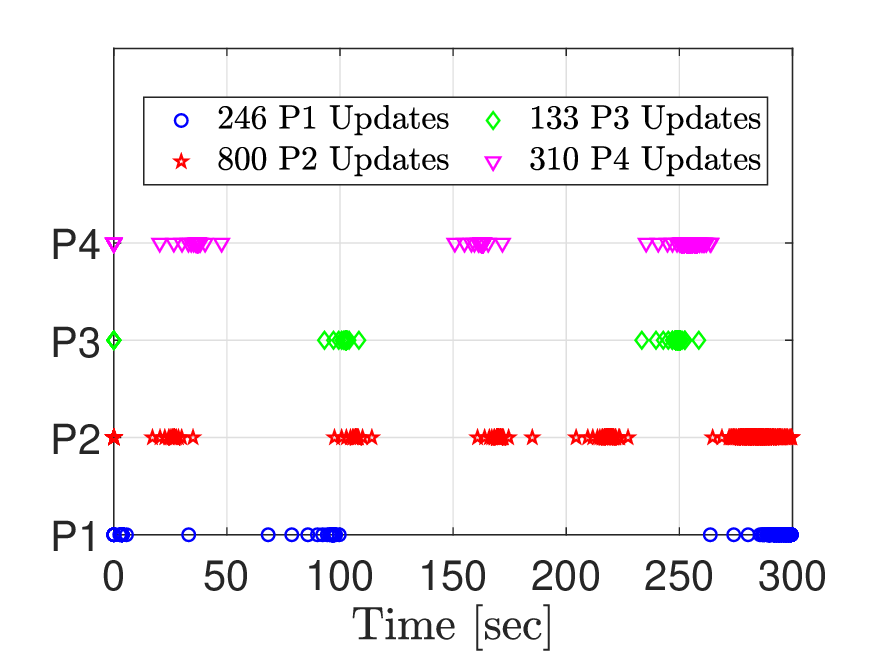}}
	\\
	\subfigure[Input signals of payoff functions, $\theta(t)$ . \label{fig:theta}]{\includegraphics[width=6.27cm]{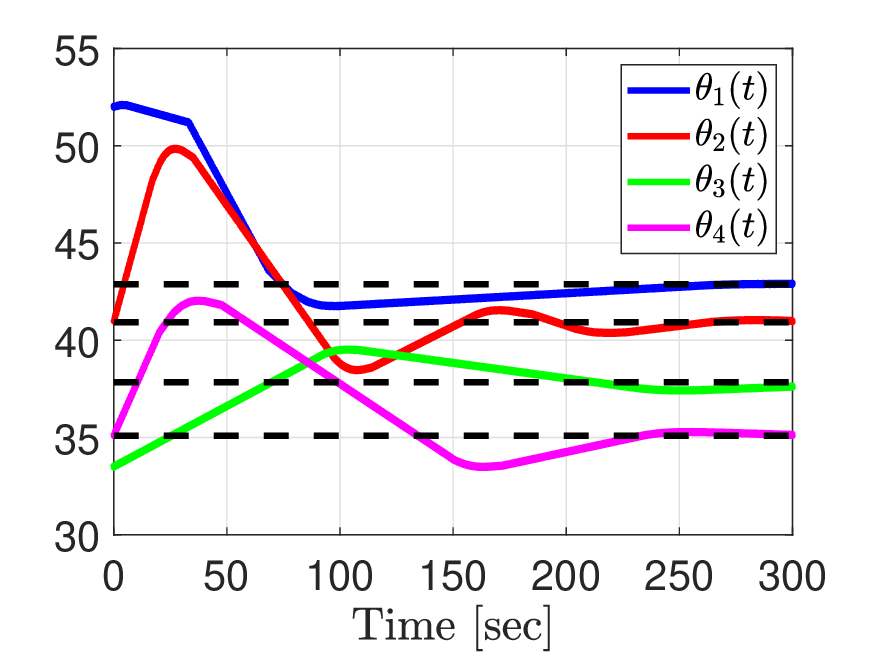}}
	~
	\subfigure[Payoff functions. \label{fig:J}]{\includegraphics[width=6.27cm]{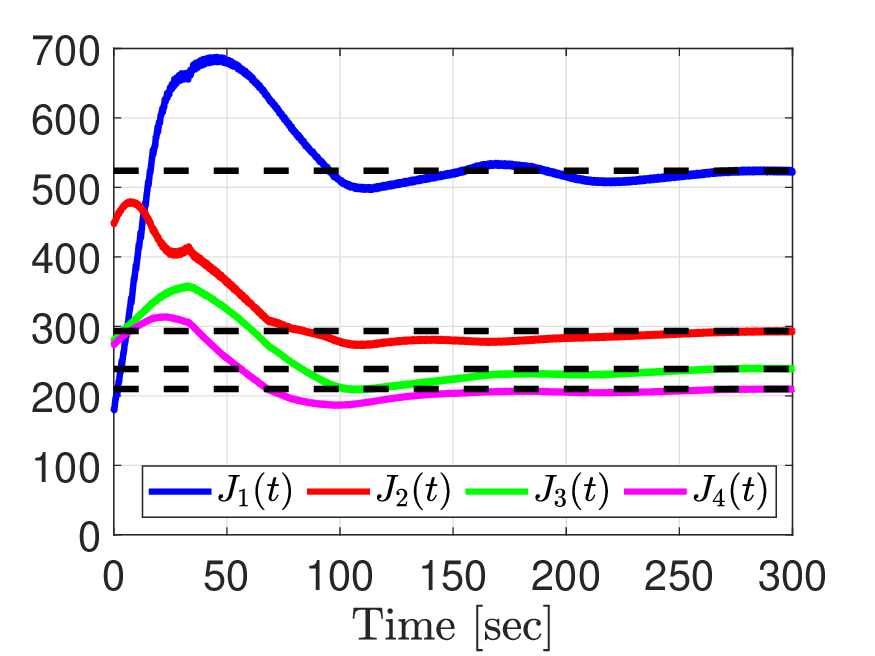}}
	\caption{Event-triggered Nash equilibrium seeking system. }\label{fig:ET_NES}
\end{figure*}

In order to achieve the Nash equilibrium (\ref{eq:theta_NE}) and (\ref{eq:J_NE}), without requiring detailed modeling information, players P1, P2, P3, and P4 implement the proposed decentralized event-triggered NES strategy to determine their optimal actions. Fig. \ref{fig:ET_NES} illustrates the corresponding time-evolution of the closed-loop system. 

Fig. \ref{fig:u} demonstrates the aperiodic piecewise-constant behavior of the players' tuning laws. Recalling that each player estimates a distinct pseudo-gradient 
 component and decides independently  when to trigger the corresponding player action, these updates occur autonomously. For the purpose of comparison, the instances of these updates are depicted in Fig. \ref{fig:controlUpdate}. In a simulation spanning 300 seconds, the actions of player P1 were updated 246 times, and those of players P2, P3, and P4 updated 
 800, 133, and 310 times, respectively. 
 
Notably, while these updates are independent, they collectively drive the
system towards the Nash equilibrium, as depicted in Fig.~\ref{fig:theta}. Moreover, Fig. \ref{fig:J} illustrates how, within the described oligopoly market structure without information sharing, each player can maximize their profits $J_{i}(t)$ by employing the proposed decentralized NES strategy via ETC to set the product prices $\theta_i(t)$ in the noncooperative game.

\newpage
\section{Conclusions}

This paper has introduced a new method for achieving locally stable convergence to Nash equilibrium in noncooperative games by employing distributed event-triggered tuning policies. This work marks the first instance of addressing noncooperative games in a model-free fashion by integrating event-triggered and extremum seeking in order to optimize the pseudo-gradient of unknown quadratic payoff functions, without sharing information. Each player evaluates independently the deviation between her own current state variable and its last broadcasted value to update her actions. Closed-loop stability and performance are preserved within the constraints of limited actuation bandwidth. The stability analysis is rigorously conducted using time-scaling techniques, Lyapunov's direct method, and averaging theory for systems with discontinuous right-hand sides. Furthermore, the paper quantifies the size of the ultimate small residual sets around the Nash equilibrium.\linebreak Simulation results considering an oligopoly market structure validate the effectiveness of the proposed approach. Future research
lies in the design and analysis of different control problems
with ET-NES, as considered in the following references \cite{OCCPBS:2017,OCH:2017,OPH:2015,POH:2019}.

\begin{ack}                               
V. H. P. Rodrigues and T. R. Oliveira would like to thank the Brazilian funding agencies CAPES, CNPq and FAPERJ for the financial support. M. Krsti{\' c} was funded under the NSF grants $CMMI-2228T91$ and $ECCS-2210315$. 
\end{ack}

\appendix

\begin{figure*}
    \centering
    \begin{adjustbox}{rotate=0, center}
        \begin{minipage}{1\linewidth}
            \centering
        \includegraphics[width=1\linewidth]{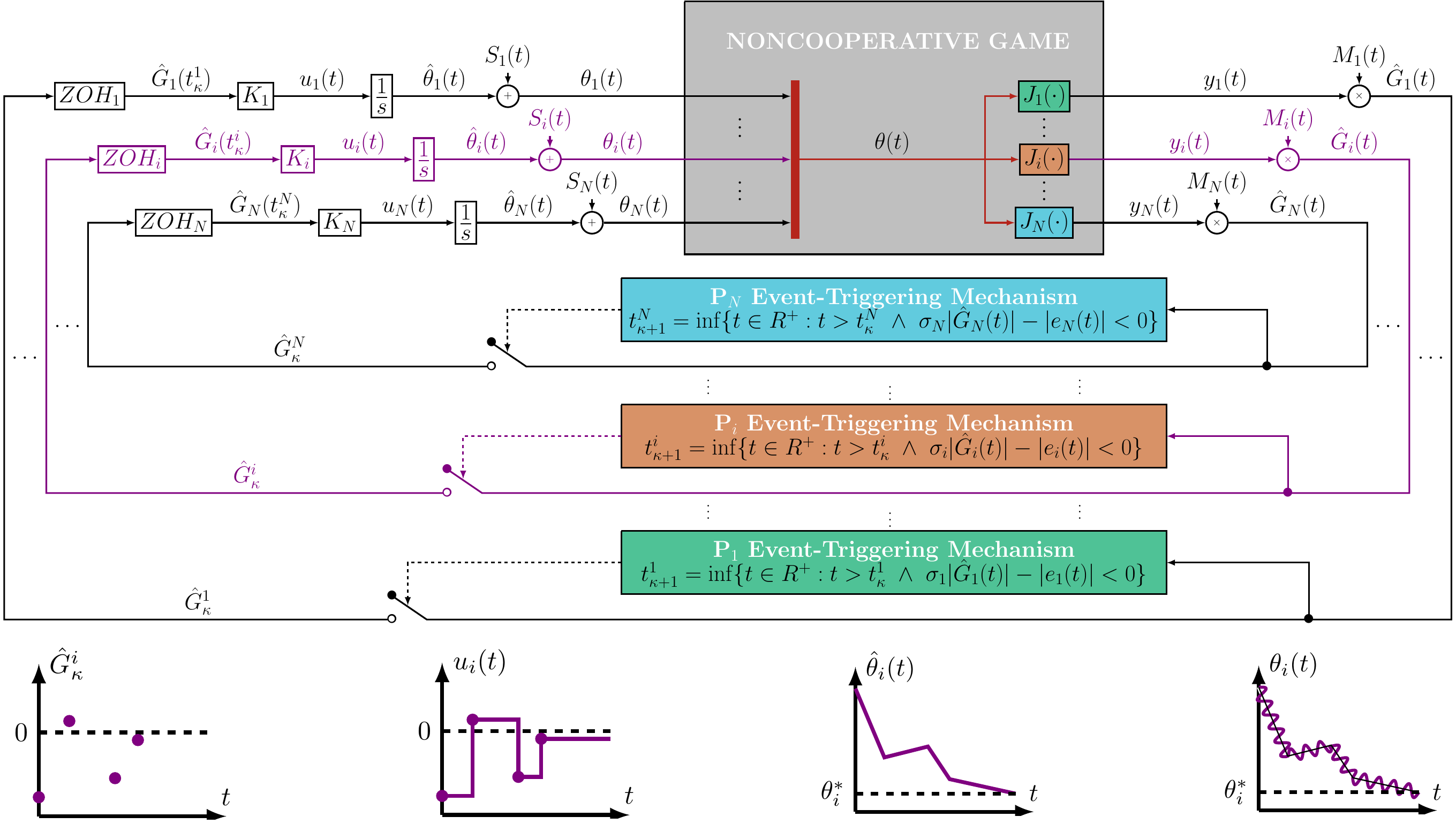}
            \label{fig:blockDiagram_2}
        \end{minipage}
    \end{adjustbox}
    \caption{Block diagram of the complete closed-loop system.}
\end{figure*}

\begin{small}
\section*{Appendix}

\section{Average of Discontinuous Systems}
\label{appendix_plotnikov}

It is important to clarify that
the following averaging result for differential inclusions with discontinuous right-hand sides by Plotnikov \cite{P:1979} takes into account discontinuities not in the periodic perturbations, but in the right-hand side of the state equation. 
In our case, perturbations of the classical ES scheme remains periodic while the discontinuities induced by the increasing sequence of event times are not periodic. As illustrated in Fig. \ref{fig:blockDiagram_2}$.1$, there is no jump in the system solutions of the closed-loop feedback. The discontinuity occurs in the 
right-hand side of the state dynamics---in the sampled-player actions $u_i(t)$. The state is continuous in such a way that it is continuously monitored by the activation mechanism. The input of the integrators in Fig. \ref{fig:blockDiagram_2}$.1$ is piecewise continuous, depending on the state. Hence, the following theorems can indeed be applied.

From \cite{P:1979}, let us consider the differential inclusion
\begin{align}
\frac{d x}{dt} \in \varepsilon X(t,x)\,, \quad x(0)=x_{0}\,, \label{eq:A1}
\end{align}
where $x$ is an n-dimensional vector, $t$ is time, $\varepsilon$ is a small parameter, and $X(t,x)$ is a multi-valued function that is $T$-periodic in $t$ and puts in correspondence with each point $(t,x)$ of a certain domain of the ($n+1$)-dimensional space a compact set $X(t,x)$ of the $n$-dimensional space. 
Let us put in correspondence with the inclusion  (\ref{eq:A1}) the average inclusion
\begin{align}
\frac{d \xi}{dt} \in \varepsilon \bar{X}(\xi)\,, \quad \xi(0)=x_{0}\,, \label{eq:A2}
\end{align}
where $\bar{X}(\xi)=\frac{1}{T}\int_{0}^{T}X(\tau,\xi)d\tau$. 

\begin{theorem} \label{thm:A1}
Let a multi-valued mapping $X(t,x)$ be defined in the domain $Q\left\{t\geq 0\,, x\in D\subset\mathbb{R}^{n}\right\}$ and let in this domain the set $X(t,x)$ be a nonempty compactum for all admissible values of the arguments and the mapping $X(t,x)$ be continuous and uniformly bounded and satisfy the Lipschitz condition with respect to $x$ with a constant $\lambda$, {\it i.e.}, $X(t,x) \subset S_{M}(0)$, $\delta(X(t,x')-X(t,x''))\leq \lambda \|x'-x''\|$, where $\delta(P,Q)$ is the Hausdorff distance between the sets $P$ and $Q$, {\it i.e.}, $\delta(P,Q)=\min\left\{d|P\subset S_{d}(Q), Q \subset S_{d}(P)\right\}$, $S_{d}(N)$ being the d-neighborhood of a set $N$ in the space $\mathbb{R}^{n}$; the mapping $X(t,x)$ be $T$-periodic in $t$; for all $x_{0}\in D'\subset D$ the solutions of inclusion (\ref{eq:A2}) lie in  the domain $D$ together with a certain $\rho$-neighborhood. Then, for each $L>0$ there exist $\varepsilon^{0}(L)>0$ and $c(L)>0$ such that for $\varepsilon \in (0\,,\varepsilon^{0}\rbrack$ and $t \in \lbrack 0, L\varepsilon^{-1}\rbrack$:
\begin{enumerate}
\item for each solution $x(t)$ of the inclusion (\ref{eq:A1}) there exists a solution $\xi(t)$ of the inclusion (\ref{eq:A2}) such that 
\begin{align}
\|x(t)-\xi(t)\|\leq c\varepsilon =\mathcal{O}(\varepsilon); \label{eq:A4}
\end{align}
\item for each solution $\xi(t)$ of the inclusion (\ref{eq:A2}) there exists a solution $x(t)$ of the inclusion (\ref{eq:A1}) such that  the inequality (\ref{eq:A4}) holds.
\end{enumerate}
Thus, the following estimate is valid: $\delta(\bar{R}(t),R'(t))\!\leq\! c\varepsilon \!=\!\mathcal{O}(\varepsilon)$, 
where $\bar{R}(t)$ is a section of the family of solutions of the inclusion (\ref{eq:A2}) and $R'(t)$ is the closure of the section $R(t)$ of the family of solutions of the inclusion (\ref{eq:A1}).
\end{theorem}

\begin{theorem} \label{thm:A2}
Let all the conditions of Theorem~\ref{thm:A1} and also the following condition be fulfilled: the $R$-solution $\bar{R}(t)$ of inclusion (\ref{eq:A2}) is uniformly asymptotically stable. Then there exist $\varepsilon^{0}>0$ and $c>0$ such that for $0< \varepsilon \leq \varepsilon^{0}$
\begin{align}
\delta(\bar{R}(t),R'(t))\leq c\varepsilon = \mathcal{O}(\varepsilon)\,, \quad \forall t\geq 0 \,.\label{eq:A6}
\end{align}
\end{theorem}
\end{small}


\end{document}